\RequirePackage{fix-cm}
\documentclass{svjour3}                     
\smartqed  
\usepackage{mathptmx}      
\usepackage{amsmath,amssymb,amsfonts}%
\usepackage{graphicx}
\usepackage{subfigure}
\usepackage{stmaryrd}
\usepackage[T1]{fontenc}
\usepackage{latexsym,mathrsfs}
\usepackage{float}
\usepackage{xcolor}
\usepackage{algorithm, algorithmic}
\usepackage{array,bm}
\usepackage[colorlinks,linkcolor=blue, anchorcolor=blue, citecolor=blue]{hyperref}
\usepackage{color,caption,cases}
 \usepackage{tabularx,booktabs }
\usepackage{indentfirst}
\usepackage{appendix}
\usepackage[numbers,sort&compress]{natbib}
\numberwithin{equation}{section}

\newcommand{\malA}{\mathcal{A}}
\newcommand{\malB}{\mathcal{B}}
\newcommand{\malD}{\mathcal{D}}
\newcommand{\malE}{\mathcal{E}}
\newcommand{\malF}{\mathcal{F}}
\newcommand{\malH}{\mathcal{H}}

\newcommand{\malO}{\mathcal{O}}
\newcommand{\mcrS}{\mathscr{S}}
\newcommand{\mcrL}{\mathscr{L}}
\newcommand{\mcrN}{\mathscr{N}}
\journalname{Arxiv }

\begin{document}
	
\title{ Strang splitting structure-preserving high-order compact difference schemes for nonlinear convection diffusion equations}	
\titlerunning{ Splitting structure-preserving HOC difference schemes }        
			
\author{
	Baolin Kuang    \and
	Shusen Xie \and 
	Hongfei Fu$^*$ 
}
\authorrunning{ Baolin Kuang \and Shusen Xie \and Hongfei Fu} 

\institute{ 
	Hongfei Fu		\at
	Corresponding author. School of Mathematical Sciences \& Laboratory of Marine Mathematics, Ocean University of China, Qingdao  266100, China\\
	\email{fhf@ouc.edu.cn} 
	\and
	Baolin Kuang  \at
	School of Mathematical Sciences, Ocean University of China, Qingdao  266100, China\\
	\email{blkuang@stu.ouc.edu.cn} 
	\and 
	Shusen Xie \at
	School of Mathematical Sciences \& Laboratory of Marine Mathematics, Ocean University of China, Qingdao  266100, China\\
	\email{shusenxie@ouc.edu.cn}
}

\date{Received: date / Accepted: date}

\maketitle

\begin{abstract}
In this paper, we present a class of high-order and efficient compact difference schemes for nonlinear convection diffusion equations, which can preserve both bounds and mass. For the one-dimensional problem, we first introduce a high-order compact Strang splitting scheme (denoted as HOC-Splitting), which is fourth-order accurate in space and second-order accurate in time. Then, by incorporating the Lagrange multiplier approach with the HOC-Splitting scheme, we construct two additional bound-preserving or/and mass-conservative HOC-Splitting schemes that do not require excessive computational cost and can automatically ensure the uniform bounds of the numerical solution. These schemes combined with an alternating direction implicit (ADI) method are generalized to the two-dimensional problem, which further enhance the computational efficiency for large-scale modeling and simulation. Besides, we present an optimal-order error estimate for the bound-preserving ADI scheme in the discrete $L_2$ norm.  Finally, ample numerical examples are presented to verify the theoretical results and demonstrate the accuracy, efficiency, and effectiveness in preserving bounds or/and mass of the proposed schemes.


\end{abstract}

\section{ Introduction}\label{sec:Intro}
 Convection diffusion  equations arise  in many science and engineering applications, such as  the transport of groundwater in soil  \cite{assouline2013infiltration},
the phenomenon of chemotaxis in biology \cite{lu2023finite},
the dispersion of tracers in anisotropic porous media \cite{fattah1985dispersion}, two-phase flow problems \cite{chen2024second} and so on.  In this paper, we are interested in efficient and reliable modeling of the following nonlinear convection diffusion equation
\begin{align}\label{model:e1}
	\left\{     
	\begin{aligned}
		& u_t +  \nabla \cdot \bm F( u) = \gamma \Delta   u +  S(\bm x,t),   &
		\quad \textrm{in} \quad \Omega \times  (0, T],  \\
		&  u(\bm x, 0) =  u^o(\bm x),  & \quad \textrm{in} \quad   \Omega \times  (0, T],  
	\end{aligned}
	\right.
\end{align}
where $\Omega=[0,L]^d \subset \mathbb{R}^d$ ($d=1,2$) is a bounded rectangular domain with $\partial\Omega$ its boundary.  The symbols $\Delta$ and $\nabla$ denote the Laplace and gradient operators, respectively. For simplicity of presentation, we assume that model \eqref{model:e1} is enclosed with periodic boundary conditions. 
Assume that $\bm F( u)$ and $S(\bm x,t)$ are well-defined smooth functions, and $\gamma>0$ denotes the diffusive coefficient. Under periodic boundary conditions, there holds 
\begin{align}  \label{model:mass}
	\int_\Omega   u(\bm x, t) d \bm x = \int_\Omega  u^o(\bm x) d \bm x+ \int_0^t \int_\Omega  S(\bm x,s)  d \bm x ds,  
\end{align}
which means \eqref{model:e1} itself is mass conservative. Therefore, this requires that our numerical schemes can also conserve the mass in the discrete level.

When diffusion  dominates  convection, standard numerical methods usually work quite satisfactorily. However, in many applications, $\gamma$ is relatively small compared to $|\bm F^{\prime}|$, so the equation becomes convection dominant. In such case, the model often exhibits strong hyperbolic characteristics, and standard methods will produce non-physical oscillation or dispersion solutions  and bring severe restriction on the temporal-spatial step-ratio. Thus, it leads to huge challenge in numerical simulations. To tackle this issue,  a variety of numerical schemes and numerical techniques have been designed and introduced to give a better simulation of such kinds of models, see \cite{shih1999modified,rui2010mass,cavalli2009family} and references therein for details. High-order compact (HOC) difference method is one of popularly used approaches to achieve high order spatial accuracy, which features smaller stencils (therefore, requiring less memory) and better resolution for high frequency waves \cite{lele1992compact}. 
Some HOC and  HOC alternating direction implicit (HOC-ADI) difference methods \cite{tian2007high,berike2007,hu2023efficient,you2006high,karaa2004high} have been proposed for linear convection diffusion equations. However, the extension of the above-mentioned high-order schemes to nonlinear unsteady convection diffusion equations is not so effective and efficient. That is one of the main objectives of this paper.

In many practical applications, solutions are always nonnegative or bounded, such as density and probability distribution. Thus, it is natural to design positivity/bound-preserving numerical schemes that maintain the interested physical quantity $u$ positive or within the interval $[m, M]$, where $m=\min ~ u^{o}$ and $M=\max ~ u^{o}$. Otherwise, the solutions are non-physical and meaningless, and may lead to instability for long-term modeling. In recent years, substantial efforts have been devoted to the construction of high-order maximum-principle-preserving schemes for convection diffusion equations \cite{xiong2015high,yang2016high,li2018high}, compressible miscible displacements \cite{guo2020high}, incompressible wormhole propagation \cite{liu2021high}, gas flow in porous media \cite{kou2023efficient}, and so on. As mentioned earlier, mass conservation \eqref{model:mass} is another important property of convection diffusion equations, and
in \cite{rui2010mass,fu2012mass,fu2017time,futai2022mass,colera2020nearly}, some mass-conservative characteristic finite element/finite volume methods were  proposed and applied for linear convection diffusion equations. Up to now, there are relatively few papers that focus on both bound-preserving and mass-conservative schemes \cite{xiong2015high,li2018high,qin2023positivity}, of which \cite{li2018high}  discussed the weak monotonicity of the HOC difference scheme  for the first time, and a simple limiter was introduced to enforce the bound-preserving property without losing mass conservation and high-order accuracy.  However, the proposed scheme is explicit and hence brings strict step-ratio restriction due to stability. This inspires us to develop some efficient semi-implicit bound-preserving and mass-conservative HOC schemes for the nonlinear model \eqref{model:e1} to further weaken the restriction. 

Recently,  Xu et al. \cite{van2019positivity} proposed a positivity-preserving high-order discontinuous Galerkin algorithm with Karush-Kuhn-Tucker (KKT)-limiter,  and a semi-smooth Newton method was employed for the optimization during time marching. To avoid solving such constrained minimization problem, by adopting the Lagrange multiplier approach and predictor-corrector method, Shen  et al. \cite{cheng2022new1,cheng2022new2} constructed some positivity/bound-preserving and mass-conservative schemes for semilinear and quasi-linear parabolic equations, where in the prediction step,  a generic semi-implicit numerical scheme is solved,  and then, structure-preserving is enforced via  Lagrange multipliers with negligible computational cost in the correction step. However, its error analysis for the Lagrange multiplier approach is challenging, and yet, few work have focused on convergence analysis of bound-preserving and mass conservation numerical schemes in this direction \cite{cheng2022new2,tong2024positivity}.
As is well known, Strang splitting method \cite{strang1968construction} is one of the most popular operator splitting methods, which has been proven to be advantageous in numerically modeling various physical problems due to its simplicity and efficiency, see \cite{li2022stability,gradinaru2008strang,liu2022second,Chertock2010p} for reference. However, to the best of our knowledge, there is yet no Strang splitting scheme which concerns both structure-preserving property and high-order accuracy for nonlinear convection diffusion equations. Inspired by \cite{cheng2022new1,cheng2022new2}, we shall focus on the combinations of the Lagrange multiplier approach and the Strang splitting method to construct some high-order structure-preserving schemes for the nonlinear convection diffusion problems. In particular, it is highly nontrivial to establish error estimates for the structure-preserving splitting scheme, due to additional difficulty arising from the splitting processes.
 
The main idea of our method is to first split the model \eqref{model:e1} into two subequations: the nonlinear convection equation and the linear diffusion equation, and then treat them separately via different temporal and spatial discretization methods.  For the former equation, we use the second-order explicit strong-stability-preserving Runge-Kutta (SSP-RK2) method \cite{shu1988efficient} for temporal discretization, which is stable under a mild CFL condition. While for the latter equation, we employ the implicit Crank-Nicolson scheme for temporal discretization and use ADI approach for multidimensional problem to enhance the overall computational efficiency. For spatial discretization, fourth-order HOC difference approximations are used for the first- and second-order spatial derivatives. Meanwhile, the Lagrange multiplier approach and predictor-corrector method are adopted to ensure the inherent bound-preserving or/and mass-conservative properties of the model. In summary, we propose some efficient and structure-preserving HOC difference schemes for the nonlinear convection diffusion equation, and our new schemes enjoy the following remarkable advantages: 
\begin{itemize}
   \item They are bound-preserving or/and mass-conservative, and automatically ensure the uniform bounds of the numerical solution;  
   
   \item The resulting algebraic systems  are linear and with  only  (cyclic)  tridiagonal constant-coefficient matrices, which largely and efficiently reduce the computational cost. Besides, the TVB limiter is also utilized to  avoid spurious  oscillations  for the transport equation. Moreover, the substeps method can further weaken the step-ratio restriction caused by the explicit discretization of the nonlinear convection, and thus further improves the computational efficiency;

  \item The ADI schemes for the multidimensional problem further enhance the overall computational efficiency; 

     \item The BP-HOC-ADI-Splitting scheme is proven to be fourth-order accurate in space and second-order accurate in time under the discrete $L_2$ norm, see Sect. \ref{sec:ErrEst}.
\end{itemize}

The remainder of the paper is organized as follows. In Sect. \ref{sec:1D}, the HOC-Splitting scheme is firstly presented for the one-dimensional (1D) problem and the discrete mass-conservative property is proved. Then, based on the proposed HOC-Splitting scheme and by adopting the Lagrange multiplier approach, the bound-preserving or/and mass-conservative HOC-Splitting schemes, denoted as BP-HOC-Splitting and BP-MC-HOC-Splitting, are further developed. In Sect. \ref{sec:2D:ADI}, some efficient ADI splitting schemes are developed for the two-dimensional (2D) model problem, and we carry out the discrete $L_2$-norm error estimates for the BP-HOC-ADI-Splitting scheme in Sect. \ref{sec:ErrEst}. Several numerical experiments are presented to illustrate the accuracy and effectiveness of the proposed schemes in Sect. \ref{sec:test}, and meanwhile, comparisons of computational efficiency as well as the structure-preserving property with other schemes are also tested. Concluding remarks are given in the last section.

\section{ High-order splitting schemes for 1D model}\label{sec:1D}
In one-dimensional case, model \eqref{model:e1} reduces to the following time-dependent nonlinear scalar convection diffusion equation with periodic boundary condition 
\begin{align}\label{mod1_1}  
	\left\{      
	\begin{aligned}
		& u_t + f(u)_x - \gamma  u_{xx}  = S(x,t),  & \quad   0<  x < L, ~ 0<   t \leq T,   \\	
		& u(x, 0)  =u^{o}  (x),  & \quad		   0\leq  x \leq L,\\
		& u(0, t)  =u(L, t), & \quad 0\leq   t \leq T.
	\end{aligned}
	\right.
\end{align}
		
In this section, we focus on constructing some high-order and efficient HOC-Splitting schemes for \eqref{mod1_1}.
  	
	\subsection{ Operator splitting and temporal discretizations}\label{sec_1d_time}
    
For a positive integer $N_t$, we introduce a uniform temporal partition of $[0,T]$ as $t_n = n \tau$ $(n=0,1,\ldots,N_t)$ with time stepsize $\tau = {T}/{N_t}$. Let $t_{n+1/2} = (t_n + t_{n+1})/2$. Moreover, we denote the solution operator associated with the diffusion process 
	\begin{align}
		{u}_t= \gamma u_{xx}  \label{1d_para}
	\end{align}
of equation \eqref{mod1_1} by $\mcrS_{\mcrL}$, and the solution operator associated with the nonlinear transport process
	\begin{align}
		{u}_t+ f(u)_x = S(x,t) \label{1d_hyper}
	\end{align}
of equation \eqref{mod1_1} by $\mcrS_{\mcrN}$.	Then, the second-order Strang splitting of  \eqref{mod1_1} involves three substeps for approximating the solution at time $t_{n+1}=t_n+\tau$  from input solution $u(x,t_n)$ \cite{strang1968construction}:
	 \begin{align}  
  {u}( {x}, t_n+\tau)=\mcrS_{\mcrL}(\tau / 2) \mcrS_{\mcrN}(\tau) \mcrS_{\mcrL}(\tau / 2)\,  {u}( {x}, t_n)  + \malO (\tau^3). \label{split_1}
	\end{align}
That is,  the second-order splitting approximation of  model \eqref{mod1_1} is reduced to the following process:
     given $u^0(x) = u^o(x)$, for $n = 0,1,\cdots, N_t-1$, find $u^{n+1}(x) = u^{***}(x,t_{n+1})$ such that
\begin{numcases}{}
  u^{*}_{t} -\gamma u^{*}_{xx} = 0,  & $t\in [t_n, t_{n+1/2}]; \quad u^{*}({x,t_n})=u^n(x)$,\label{strang_1} \\
		u^{**}_{t} + f(u^{**})_x = S(x, t),   & $t\in [t_n, t_{n+1}]; \qquad u^{**}(x,t_n)=u^{*}(x,t_{n+{1}/{2}})$, \label{strang_2} \\
		u^{***}_{t} -\gamma u^{***}_{xx} = 0,  &  $t\in [t_{n+1/2}, t_{n+1}]; \quad u^{***}(x,t_{n+{1}/{2}})=u^{**}(x,t_{n+1})$. \label{strang_3}
\end{numcases}	
\begin{remark} 
Of course, at each time interval $[t_n, t_{n+1}]$, we can first calculate the transport in the first half-step $[t_n, t_{n+1/2}]$, then calculate the diffusion in the whole step $[t_n, t_{n+1}]$, and finally calculate the transport again in the second half-step $[t_{n+1/2}, t_{n+1}]$, i.e.,
		\begin{align}
			{u}( {x}, t_n+\tau)=\mcrS_{\mcrN}(\tau / 2) \mcrS_{\mcrL}(\tau) \mcrS_{\mcrN}(\tau / 2) \, u(x, t_n) + \malO (\tau^3). \label{split_2}
		\end{align}
The main differences between the two splitting methods are the stability and computational efficiency \cite{zhang2016operator}. The former splitting \eqref{split_1}, which implicitly treats the stiff diffusion term twice, usually has better stability than the latter splitting \eqref{split_2}. In contrast, the CPU time cost by the latter one is usually less than the former one. In this paper, we focus on the splitting method \eqref{split_1}.
\end{remark}		
		
\begin{remark} In practice, the exact solution operators $\mcrS_{\mcrL}$ and $\mcrS_{\mcrN}$ shall be approximated by some temporal and spatial numerical discretizations. Since the diffusion problem \eqref{1d_para} and the transport problem \eqref{1d_hyper} are of different natures, they shall be solved via different discretization methods, which is one of the main advantages of operator splitting technique. We shall adopt implicit temporal discretization for the diffusion equation due to stiffness, while for the transport equation we shall adopt explicit strong-stability-preserving temporal discretization.
\end{remark}

Next, as the splitting error of the Strang splitting method is of second-order accuracy, 
we introduce the temporal second-order semi-discretization scheme for the splitting method \eqref{strang_1}--\eqref{strang_3} as follows:
\begin{align}\label{Split:semi:e1}
	\left\{     
	 \begin{aligned}
	& \dfrac{U^{n,1}-U^{n}}{{\tau}/{2}} - \gamma \dfrac{U_{xx}^{n,1} + U_{xx}^{n} } {2} =0, & (\textrm{CN})\\
	&  \dfrac{U^{n,2} -U^{n,1} } {\tau}  +  f(U^{n,1}  ) _x = S(x,t^n), & (\textrm{SSP-RK2: P})\\ 
    & \dfrac{{U}^{n,3} -U^{n,1}}{\tau} + \dfrac{1}{2}     \big(f(U^{n,1}) _x+f(U^{n,2}  ) _x\big) 
	= \dfrac{1}  {2}  \big( S(x,t^n) + S(x,t_{n+1})   \big),  & (\textrm{SSP-RK2: C}) \\    
    & \dfrac{U^{n+1}-U^{n,3}}{ {\tau}/{2}} -\gamma \dfrac{U_{xx}^{n+1} +{U}_{xx} ^{n,3}}{2}  = 0, & (\textrm{CN})
	 \end{aligned}
	 \right.
	\end{align}
where, the implicit Crank-Nicolson (CN) discretization is adopted for the diffusion process to obtain the intermediate variable $U^{n,1}$ and the solution $U^{n+1}$, and the explicit two-stage SSP-RK2 discretization \cite{shu1988efficient} is adopted for the nonlinear transport process to obtain the intermediate variables $U^{n,2}$ and $U^{n,3}$. The explicit temporal discretization, which consists of two simple predictor (P) and corrector (C) stages, indeed simplifies the computation for the nonlinear transport process \eqref {strang_2}, and of course it also inevitably brings stability restriction for the time stepsize, but this is acceptable compared with direct explicit discretization of model \eqref{mod1_1} \cite{gottlieb2001strong}.

\subsection{ HOC spatial discretization}\label{sec_1d_space}
Given a positive integer $N_x$, denote  the sets of spatial grids $\Omega_h=\left\{ x_i = i h_x  \mid 0 \leq i \leq N_x \right\}$ with stepsize $h_{x} = {L} / {N_x}$. Accordingly, define the space of periodic grid function on $\Omega_h$ by $\mathcal{V}_h = \left\{u= \{u_{i}\}\mid u_{i+N_x} = u_i \right\} $. For any grid function $v \in \mathcal{V}_h $, we introduce the standard difference operators and fourth-order compact operators as follows:
		\begin{equation*}
			\begin{aligned}
				\malD_{\hat{x}} v_i =\frac{v_{i+1}-v_{i-1}}{2 h_x},  ~
				\delta_x^2 v_i =\frac{v_{i+1}-2 v_i + v_{i-1}}{h_x^2}, ~
    		\malA_{x} v _{i}= 
				\Big(I+\frac{h_x^2}{12} \delta_x^2 \Big) v_{i}, ~
				\malB_{x} v  _{i }=
				\Big( I+\frac{h_x^2}{6} \delta_x^2 \Big) v_{i}.
			\end{aligned}
		\end{equation*}
      It is shown in  \cite{tian2007high,berike2007,hu2023efficient,li2018high}  that
		\begin{equation}  \label{compact:AB}
			v_{x x}  = ( \malA_{x} )^{-1} \delta_x^2 v+\malO (h_x^{4}  ),  \quad	
			v_{x}   = ( \malB_{x})^{-1}\malD_{\hat{x}} v+\malO (h_x^{4} ).
		\end{equation}
     We also introduce  the discrete $L_2$ inner product $(v,w):= \sum_{i=1}^{N_x} h_x  v_i w_i$ for $v,w \in \mathcal{V}_h $.

Let $S^n_i =S( x_i, t^n )$, and we use $u_h^n=\{u_i^n\}\in \mathcal{V}_h$ to denote the fully-discrete numerical solution at $t=t^n$, and $u_h^{n,k}=\{u_i^{n,k}\}$, $1\le k \le 3$, to represent the fully-discrete intermediate solutions. By applying the compact spatial approximations in \eqref{compact:AB}, a fully-discrete HOC-Splitting difference scheme for \eqref{Split:semi:e1} is proposed as follows:
\begin{numcases}  {}  
			\label{1d:scheme:diffu1}  
			\malA_{x} \frac{u_i^{n,1 } -u_i^{n}} { {\tau} / {2}} -\gamma \delta_x^2 \frac{u_i^{n ,1 }   + u_i^{n}}{2}  =0, \\
			\label{1d:scheme:convec1}
			\malB_{x}  \frac{u_i^{n  ,2}  - u_i^{n,1  }  }  {\tau}  + \mathcal  D_{\hat{x}}  f(u_i^{n,1 }  )= \malB_{x}    S_i^{n},    \\
			\label{1d:scheme:convec2} 
			\malB_{x}   \frac{ {u}_i  ^{n,3}  -u_i^{n ,1 }  }  {\tau}  +\frac{ 1 }  {2} \malD_{\hat{x}} \big( f(u_i^{n ,1 }  )+f(u_i^{n  ,2}  ) \big)= \frac{1}  {2} \malB_{x}   \big(S_i^{n} + S_i^{n+1} \big),  \\
			\label{1d:scheme:diffu2}
			\malA_{x}   \frac{u_i^{n+1}  -  {u}_i  ^{n,3}  } { {\tau} / {2}}  -\gamma  \delta_x^2   \frac{u_i^{n+1}  + {u}_i  ^{n,3}  }  {2}  =0.
\end{numcases}

\begin{remark}\label{rem:e1}  The proposed HOC-Splitting scheme \eqref{1d:scheme:diffu1}--\eqref{1d:scheme:diffu2} exhibits second-order accuracy in time and fourth-order accuracy in space.  Note that \eqref{1d:scheme:convec1} and \eqref{1d:scheme:convec2} are no longer explicit compared to their semi-discrete form \eqref{Split:semi:e1}, but they can be solved very easily and fast. In fact, the coefficient matrices of the resulting linear algebraic systems corresponding to \eqref{1d:scheme:convec1} and \eqref{1d:scheme:convec2} are symmetric positive definite and with the same cyclic tridiagonal structure generated by the compact operator $\malB_{x}$. On the other hand, \eqref{1d:scheme:diffu1} and \eqref{1d:scheme:diffu2} own the same cyclic tridiagonal matrix generated by the operator $\big(\malA_{x}-\frac{\gamma \tau}{4} \delta_x^2\big)$. As a consequence, both the memory requirement and computational complexity of the scheme are of only $\malO (N_x)$ per time step. Thus, the proposed scheme is well-suited for practical numerical simulations. However, the above scheme cannot guarantee the bound-preserving property for the physical quantity.
\end{remark}

 \begin{remark}\label{rem:1d:K}  Note that the  CFL condition only comes from the transport process  \eqref{1d:scheme:convec1}--\eqref{1d:scheme:convec2}. To  weaken the time step-ratio restriction, one can resort to the substeps method \cite{lan2023operator} to evolve the transport equation from $t_n$ to $t_{n+1}$ as needed. To be specific, we run \eqref{1d:scheme:convec1}--\eqref{1d:scheme:convec2} that originally acts on the time interval $[t_n, t_{n+1}]$ over $K$ subintervals $[t_{n,k}, t_{n,k+1}]=[t_n  + k\delta t, t_n  + (k+1)\delta t]$ $(k =0, \ldots,K-1)$ with a small time stepsize $ \delta t  = \tau/ K $  by
   	\begin{align}  \label{substep_1}
   	\left\{    
   	\begin{aligned}
  	 &  \malB_{x}  \frac{u_{i,k+1}^{n  ,2}  - u_{i,k}^{n,1  }  }  {\delta t}  + \mathcal  D_{\hat{x}}  f(u_{i,k}^{n ,1 }  )= \malB_{x}    S_{i,k}^{n},  \\
     &  \malB_{x}  \frac{ {u}_{i,k+1}  ^{n,3}  -u_{i,k}^{n ,1 }  }  {\delta t}  +\frac{ 1 }  {2} \malD_{\hat{x}} \big(  f(u_{i,k}^{n ,1 }  )+f(u_{i,k+1}^{n  ,2}  ) \big)= \frac{1}  {2} \malB_{x}   \big(S_{i,k}^{n} + S_{i,k+1}^{n} \big),  
   	\end{aligned}
   	\right.
   \end{align}
where  
$u_{i,0}^{n,1}:=u_{i}^{n,1}$ and $S_{i,k}=S(x_i, t_{n,k})$. In this context, for the implementation of the HOC-Splitting scheme,  one needs to solve $2K$ linear algebraic systems resulting from explicit treatment of the transport equation with a smaller time stepsize $ \delta t$ and two linear algebraic systems resulting from the implicit discretization of the diffusion equation with a larger time stepsize $\tau$.  Below, we abbreviate the scheme \eqref{1d:scheme:diffu1}, \eqref{substep_1}   and \eqref{1d:scheme:diffu2}  as HOC-Splitting-Substeps. We shall illustrate the effectiveness of the substeps method in Example \ref{1d:test:burger}, where numerical results show that this approach does not affect the accuracy, but can efficiently reduce the time step-ratio restriction and therefore improves the computational efficiency.
    \end{remark}
\begin{remark}\label{rem:1d:limter}  
	To achieve nonlinear stability and eliminate non-physical oscillations for the transport equation  \eqref{1d:scheme:convec1}, a TVB limiter introduced in \cite{cockburn1994nonlinearly} can be adopted.  For simplicity, we consider $S=0$ and denote $\overline u =\malB_{x} u$.  Then, \eqref{1d:scheme:convec1}  can be written into a conservation form:
		 \begin{align}\label{1d:scheme:convec1:tvb} 
		 	\overline{u}_i^{n,2} =\overline{u}_i^{n,1}
                     -\frac{\tau}{h_x} \big(\hat{f}_{i+\frac{1}{2}} -\hat{f}_{i-\frac{1}{2}} \big), \quad 
          \hat{f}_{i+\frac{1}{2}}=\frac{1}{2} \big(f(u_{i+1}^{n,1})+f(u_i^{n,1})\big).
		 \end{align}
Using the TVB limiter to get the modified flux $\hat{f}_{i\pm\frac{1}{2}}^{(m)}$ to approximate $\hat{f}_{i\pm\frac{1}{2}}$, numerical oscillations can be effectively reduced; see \cite{cockburn1994nonlinearly} for details. In the same way, the TVB limiter can also be given for \eqref{1d:scheme:convec2}. In the subsequent Example \ref{1d:test:burger}, we shall also verify the effectiveness of this method.
\end{remark}	

Finally, we end this subsection by proving the discrete version mass conservation relation \eqref{model:mass}.
\begin{theorem}(Discrete Mass Conservation) \label{thm:mass:1d}
    Let $u_h^n= \{u_{i}^{n}\} \in \mathcal{V}_h$ be the solution of the HOC-Splitting scheme \eqref{1d:scheme:diffu1}--\eqref{1d:scheme:diffu2}. Then there holds 
        \begin{align}\label{MassCons:D}
	        ( u_h^{n+1}, 1 ) =( {u}  ^{o} , 1 )+ \frac{\tau}  {2} \sum_{k=0}^{n}  (S^{k} + S^{k+1}  , 1 ). 
        \end{align}
\end{theorem}    
\begin{proof}
   First, taking the inner product on both sides of \eqref{1d:scheme:diffu1} with 1, we obtain
        \begin{align*}  
        	\Big( \malA_{x} \frac{u_h^{n,1 } -u_h^{n}} { {\tau} / {2}}, 1 \Big)  
                    - \Big( \gamma    \delta_x^2 \frac{u_h^{n ,1 }   + u_h^{n}}{2}, 1 \Big)=0. 
        \end{align*}
      As the periodicity of the solution, we see
        \begin{align*}
        \big(\delta_x^2 v,1\big)  = 0\quad \text{and} \quad
        \big( \malA_{x} v, 1\big)   =  ( v, 1) +  \frac{h_x^2}{12} \big ( \delta_x^2 v, 1\big)  =  (v,1),
        \end{align*}
        which gives
         \begin{align}\label{MassCons:D1}
        	  (u_h^{n,1}, 1) = (u_h^{n}, 1).
        \end{align}    
    Similarly,  \eqref{1d:scheme:diffu2} yields    
    \begin{align}\label{MassCons:D2}
    	 ( u_h^{n+1}, 1 ) =( {u}_h^{n,3} , 1 ).
    \end{align}        
   
   Next, taking the inner product on both sides of \eqref{1d:scheme:convec2}  with 1, and noting that 
        \begin{align*}
        (\malD_{\hat{x}} v,1)  = 0 \quad \text{and} \quad
        ( \malB_{x} v, 1 )   =   ( v, 1) +  \frac{h_x^2}{6} ( \delta_x^2 v, 1 ) =  (v,1), 
        \end{align*}
 we have
          \begin{align}\label{MassCons:D3}
        	 ( {u}_h^{n,3} , 1 )  =  (u_h^{n,1}, 1  ) + \frac{\tau}{2} (S^{n} + S^{n+1} , 1).
        \end{align}    
 
    Now, collecting \eqref{MassCons:D1}--\eqref{MassCons:D3} together, we obtain
            \begin{align}\label{MassCons:D4}
            	(u_h^{n+1}, 1 ) =( {u}_h^{n} , 1 )+ \frac{\tau}{2} (S^{n} + S^{n+1} , 1), 
           \end{align}    
      which directly implies the conclusion by adding up $n$.     \qed
\end{proof} 
       
 \subsection{Bound-preserving or/and mass-conservative HOC-Splitting schemes}\label{subsec:1d:lag}
As mentioned above, in many situations, the physical quantity $u$ of model \eqref{mod1_1} has either bound-preserving or mass-conservative property or both. However, the proposed scheme can only achieve mass conservation as proved in Theorem \ref{thm:mass:1d}, but it may not necessarily preserve the bounds. Recently, by introducing Lagrange multipliers, Shen et al. \cite{cheng2022new1,cheng2022new2} constructed some efficient positivity/bound-preserving or/and mass-conservative schemes for semilinear and quasilinear parabolic equations. In this subsection, by incorporating the proposed HOC-Splitting scheme with the Lagrange multiplier approach, we intend to construct some bound-preserving or/and mass-conservative HOC-Splitting schemes for the nonlinear convection diffusion equation \eqref{mod1_1}.
			
\subsubsection{ BP-HOC-Splitting scheme }\label{subsubsec:1d:lag_bound}
In order to motivate our new  schemes, we briefly review below the Lagrange multiplier approach proposed in \cite{cheng2022new2} combined with the HOC spatial discretization.	
To this end, we first define a quadratic function $H(u)=(M-u)(u-m)$, and it suffices to solve the following expanded system with Lagrange multiplier
	 \begin{align}  \label{semi1}
	 	\left\{     
	 	\begin{aligned}
	  & \partial_t u_i
	  + \gamma \malA_{x}^{-1} \delta_ x^2   u_i + \mathcal{B}_x^{-1} D_{\hat{x}} f(u_i) = S(x_i,t ) + \lambda_i(t) H^{\prime}(u_i), \quad  \forall u_i \in \mathcal{V}_h,  \\
	 &\lambda_i(t) \geq 0, ~ H(u_i) \geq 0, ~ \lambda_i(t) H (u_i)=0,  
	 	\end{aligned}
	 \right.
 \end{align} 
 to preserve the bounds $u_i \in [m,M]$. Equation \eqref{semi1} represents the well-known Karush-Kuhn-Tucker (KKT) condition for constrained minimization \cite{KKTBerg} and can also be viewed as a semi-discrete HOC approximation to \eqref{mod1_1}.   
 
Existing approaches for simulations of \eqref{semi1} usually start with an implicit temporal discretization scheme, so that nonlinear constrained minimization systems via iterative methods have to be solved at each time step. Instead, the authors \cite{cheng2022new1,cheng2022new2} developed an alternative efficient approach, which decouples the computation of Lagrange multiplier from that of the solution. They proposed some $k$th-order bound-preserving predictor-corrector schemes based on the backward difference formula (BDF) and Adams-Bashforth extrapolation method for temporal discretization. Such ideas can also be directly applied to model \eqref{mod1_1} combined with the HOC spatial discretization, which consists of the following two steps:
\paragraph{Step 1 (Predictor).} Solve $\tilde{u}_h^{n+1}=\{\tilde{u}_i^{n+1}\}$ via the second-order BDF2-HOC method
 \begin{equation}\label{1d:BDF2-HOC:pre}
	\begin{aligned}
		& \frac{ 3 \tilde{u}_i^{n+1}-  4u_i^n  + u_i^{n-1}  }{2\tau } - \gamma \malA_{x}^{-1} \delta_ x^2   \tilde{u}_i^{n+1} + \mathcal{B}_x^{-1} D_{\hat{x}} f(2u_i^{n}-u_i^{n-1}) =S_i^{n+1} + \lambda_i ^n H^{\prime}(u_i^{n }).
	\end{aligned}
 \end{equation}
\paragraph{Step 2 (Corrector).} Solve $u_h^{n+1}=\{u_i^{n+1}\}$ and $\lambda_h^{n+1}=\{\lambda_i^{n+1}\}$ via the  KKT  condition   
\begin{equation}\label{1d:BDF2-HOC:correct}
\left\{   
	\begin{aligned}
    & \frac{ 3 u_i^{n+1}-3\tilde{u}_i^{n+1} }{2\tau}=\lambda_i^{n+1} H^{\prime} (u_i^{n+1})  -  \lambda_i^n  H^{\prime}\left(u_i^n \right) ,  \\
	& \lambda_i^{n+1}  \geq 0, ~ H(u_i^{n+1}) \geq 0, ~\lambda_i^{n+1}  H(u_i^{n+1})=0. 
	\end{aligned}
\right.
\end{equation}
This step can also be thought of as a correction to $\tilde{u}_h^{n+1}$ to preserve the bounds $u_h^{n+1}\in [m,M]$ using the Lagrange multiplier $\lambda_h^{n+1}$, while accuracy can still be maintained. For positivity preservation, we can simply choose $H(u) = u$ in the above scheme. For simplicity, below, we abbreviate the above bound-preserving HOC scheme as BP-BDF2-HOC. 

Although the above BP-BDF2-HOC scheme can well preserve the bounds, there is still a lack of high efficiency and stability for long-term and large-scale modeling. Therefore, combined with the proposed HOC-Splitting scheme \eqref{1d:scheme:diffu1}--\eqref{1d:scheme:diffu2} with the Lagrange multiplier approach, a novel bound-preserving HOC-Splitting scheme, namely BP-HOC-Splitting, is proposed as follows: 
\paragraph{Step 1 (Predictor).}  Solve $\overline{u}_h^{n+1}$ via the HOC-Splitting method
\begin{equation} \label{scheme:1d:BP:pre}
	\left\{     
	\begin{aligned}
    &\malA_{x}   \frac{u_i^{n,1  }  -u_i^{n} }{ {\tau} / {2}  }  -\gamma\delta_x^2 \frac{u_i^{n ,1 }  +u_i^{n}  }  {2}  = \frac{1}{2}   \malA_{x} (\lambda_i^n H^{\prime}(u_i^n)),      \\
    &\malB_{x}  \frac{u_i^{n  ,2}  - u_i^{n,1  }  }  {\tau}  + \mathcal  D_{\hat{x}}  f(u_i^{n ,1 }  )= \malB_{x}    S_i^{n},   \\
    &\malB_{x}   \frac{ {u}_i  ^{n,3}  -u_i^{n ,1 }  }  {\tau}  +\frac{ 1 }  {2} \malD_{\hat{x}} \big(f(u_i^{n ,1 }  )+f(u_i^{n  ,2} )\big)= \frac{1}  {2} \malB_{x}   \big(S_i^{n} + S_i^{n+1} \big),  	\\
    &\malA_{x}   \frac{ \overline u_i^{n+1}  -  {u}_i  ^{n,3}  }  { {\tau} / {2}  }  -\gamma \delta_x^2   \frac{ \overline u_i^{n+1}  + {u}_i ^{n,3}  }  {2}  = \frac{1}{2} \malA_{x} (\lambda_i^n H^{\prime}(u_i^n)).    
	\end{aligned}
	\right.
\end{equation}
\paragraph{Step 2 (Corrector).} 	Update $u_h^{n+1}$ and $\lambda_h^{n+1}$ via the KKT condition	  		
\begin{equation}\label{scheme:1d:BP:cor}
	\left\{   
	\begin{aligned}
	    &	\frac{ u_i^{n+1} -\overline{u}_i^{n+1}  }{\tau } =  \lambda_i^{n+1} H^{\prime}(u_i^{n+1} )- \lambda_i^n H^{\prime}(u_i^n ), \\
	  & \lambda_i^{n+1}  \geq 0, ~	H(u_i^{n+1}) \geq 0, ~ \lambda_i^{n+1} H(u_i^{n+1})=0.
	\end{aligned}\right.
\end{equation}
			
In the prediction step,  the efficient Strang splitting HOC method is utilized to achieve high-order accuracy. However, the obtained solution $\overline u_i^{n+1}$ may not belong to $[m, M]$. Then, in the correction step, the Lagrange multiplier $\lambda_i^{n+1}$ is introduced to enforce the pointwise bound-preserving of $u_i^{n+1}$. In fact, the KKT condition \eqref{scheme:1d:BP:cor} can be efficiently solved as follows \cite{cheng2022new1,cheng2022new2}:
			\begin{align}\label{scheme:1d:BP:cut}
				\big(u_i^{n+1} , \lambda_i^{n+1} \big)
				=\left\{\begin{array}{cll}
					\big(\overline{u}_i^{n+1} +\eta_i^{n+1} , 0 \big), & \text { if } \quad m< \overline{u}_i^{n+1} +\eta_i^{n+1} <M,   \\[0.1in]
					\Big( m, \frac{ m-( \overline{u}_i^{n+1} +\eta_i^{n+1} ) }{\tau H^{\prime}(m)}\Big), & \text { if } \quad \overline{u}_i^{n+1} +\eta_i^{n+1}  \leq m, \\[0.1in]
					\Big(M,  \frac{ M-( \overline{u}_i^{n+1} +\eta_i^{n+1}  )}{ \tau H^{\prime}(M)}\Big), & \text { if} \quad \overline{u}_i^{n+1} +\eta_i^{n+1}  \geq M,
				\end{array}\right.  
			\end{align}
	where 
	$ \eta_i^{n+1} :=-{ \tau} \lambda_i^n H^{\prime}(u_i^n)$.

\begin{remark} \label{rm_lambda_0}
Obviously, \eqref{scheme:1d:BP:cut} implies that the above scheme is indeed bound-preserving and can also be viewed as a high-order cut-off scheme. Moreover, $\lambda_i^n$ is an approximation to zero as $h_x \rightarrow 0$, and based on the discussion in \cite{cheng2022new2}, an alternative approach is to replace $\lambda_i^n H^{\prime}(u_i^n ) $ by zero in \eqref{scheme:1d:BP:pre}--\eqref{scheme:1d:BP:cor}. In this sense, the correction step is equivalent to a very simple cut-off scheme, which is still second-order (in time) convergent with respect to the maximum norm, i.e., the temporal accuracy is not affected. Numerical experiments also verify this point.
\end{remark}
			
\subsubsection{ BP-MC-HOC-Splitting scheme} \label{subsubsec:1d:lag_mass}
A drawback of the BP-HOC-Splitting scheme \eqref{scheme:1d:BP:pre}--\eqref{scheme:1d:BP:cor}  is that it does not necessarily preserve the mass even though the original scheme \eqref{1d:scheme:diffu1}--\eqref{1d:scheme:diffu2} does. Therefore, we introduce another Lagrange multiplier $\xi^{n+1}$, which is only time-dependent, to preserve both the bounds and the mass of the physical quantity. The scheme, namely BP-MC-HOC-Splitting, is proposed as follows:		
\paragraph{Step 1 (Predictor).}  Solve $\tilde {u}_h^{n+1}$ from $u_h^{n}$ using the developed HOC-Splitting scheme 
  \begin{numcases}  {}
	\malA_{x}   \frac{u_i^{n,1  }  -u_i^{n} }{ {\tau} / {2}  }  -\gamma^x\delta_x^2 \frac{u_i^{n ,1 }  +u_i^{n}  }  {2}  = \frac{ 1 }  {2}\malA_{x}  \big(\xi^{n} + \lambda_i^n  H^{\prime}(u_i^n)\big),  \label{scheme:1d:MC:p1}   \\
	\malB_{x}  \frac{u_i^{n  ,2}  - u_i^{n,1  }  }  {\tau}  + \mathcal  D_{\hat{x}}  f(u_i^{n,1 }  )= \malB_{x}    S_i^{n},  \label{scheme:1d:MC:p2}  \\
	\malB_{x}   \frac{ {u}_i  ^{n,3}  -u_i^{n ,1 }  }  {\tau}  +\frac{ 1 }  {2} \malD_{\hat{x}}\big(f(u_i^{n ,1 }  )+f(u_i^{n  ,2}  )\big)= \frac{1}  {2} \malB_{x}   \big(S_i^{n} + S_i^{n+1} \big) , \label{scheme:1d:MC:p3}  	\\
	\malA_{x}   \frac{ \tilde u_i^{n+1}  -  {u}_i  ^{n,3}  }  { {\tau} / {2}  }  -\gamma \delta_x^2   \frac{ \tilde u_i^{n+1}  + {u}_i ^{n,3}  }  {2}  = \frac{ 1 }  {2}\malA_{x}  \big(\xi^{n} + \lambda_i^n  H^{\prime}(u_i^n)\big).   \label{scheme:1d:MC:p4}
\end{numcases}
			
\paragraph{Step 2 (Corrector).} 	  Solve $(u_h^{n+1}, \lambda_h^{n+1}, \xi^{n+1})$ via the KKT conditions			
			\begin{numcases}  {}
				\frac{u_i^{n+1}   -\tilde {u}_i^{n+1}   }{\tau } 
				=   \lambda_i^{n+1}  H^{\prime}(u_i^{n+1}   ) 
				- \lambda_i^n  H^{\prime}(u_i^n   ) +\xi^{n+1}- \xi^n, \label{scheme:1d:MC:c1}\\
				\lambda_i^{n+1}  \geq 0, ~ H(u_i^{n+1}   ) \geq 0,  ~
				\lambda_i^{n+1}  H(u_i^{n+1}   ) =0, 		\label{scheme:1d:MC:c2}\\
				(u_h^{n+1}, 1 ) = (u^n, 1 ) + \frac{\tau}  {2}  \big( S^{n} + S^{n+1} , 1\big). \label{scheme:1d:MC:c3}
				\end{numcases} 

Note that the Lagrange multiplier $\xi$ is introduced to ensure mass conservation \eqref{scheme:1d:MC:c3}. In order to solve the above system \eqref{scheme:1d:MC:c1}--\eqref{scheme:1d:MC:c3} efficiently, we denote
			\begin{align}\label{lag4}
				\theta_i^{n+1} :=  {\tau} \, \big(\xi^{n+1}- \xi^n - \lambda_i^n H^{\prime}(u_i^n) \big),
			\end{align}
			and rewrite \eqref{scheme:1d:MC:c1} as
			\begin{align*}
				\frac{u_i^{n+1}   -(\tilde{u}_i^{n+1}   +\theta_i^{n+1} )}{\tau }=  \lambda_i^{n+1} H^{\prime} (u_i^{n+1} ) .
			\end{align*}
Hence, assuming $\xi^{n+1}$ is known, we find from above and \eqref{scheme:1d:MC:c2} that
			\begin{align}
				& \big (u_i^{n+1}   , \lambda_i^{n+1}   \big)=\left\{\begin{array}{ll}
				\big(\tilde{u}_i^{n+1}   +\theta_i^{n+1}   , 0 \big), & \text { if } \quad m<\tilde{u}_i^{n+1}   +\theta_i^{n+1}   <M, \\[0.1in]
				\Big(m, \frac{  m-(\tilde{u}_i^{n+1}   +\theta_i^{n+1}  )}{\tau H^{\prime}(m)}
				\Big), & \text { if } \quad \tilde{u}_i^{n+1}   +\theta_i^{n+1}  \leq m, \\[0.1in]
				\Big(M, \frac{ M-(\tilde{u}_i^{n+1}   +\theta_i^{n+1} )   }{{\tau } H^{\prime}(M)}
				\Big), & \text { if } \quad \tilde{u}_i^{n+1}   +\theta_i^{n+1}  \geq M.
				\end{array}  \right. \label{lag5}
			\end{align}
			
As for $\xi^{n+1}$, let
    \begin{equation*}
			 { }^m \mathcal{V}_h(\theta):=\Big\{\theta_i \mid \tilde{u}_i^{n+1} + \theta_i \leq m\Big\}, ~  
                 \mathcal{V}_h^M(\theta):=\left\{\theta_i \mid \tilde{u}_i^{n+1} + \theta_i \geq M\right\}, 
    \end{equation*}
        \begin{equation*}
					 { }^m \mathcal{V}_h^M(\theta):=\left\{\theta_i \mid m<\tilde{u}_i^{n+1} +\theta_i <M\right\}.
\end{equation*}
   Then, thanks to \eqref{lag5}, the discrete mass conservation formula \eqref{scheme:1d:MC:c3} can be rewritten as 
			\begin{align}\label{lag7}
		   \sum_{ {^m\mathcal{V}_h^M}(\theta)}(\tilde{u}_i^{n+1} + \theta_i ) h_x  +\sum_{ \mathcal{V}_h^M(\theta)} M h_x  +\sum_{ {^m\mathcal{V}_h (\theta)}} m h_x = (u^n,1)  + \frac{\tau}  {2}  \big( S^{n} + S^{n+1} , 1\big).  
			\end{align}
	Setting
			\begin{align*}
				\begin{aligned}
					& G_n(\theta) :=\sum_{ {^m\mathcal{V}_h^M}(\theta)}(\tilde{u}_i^{n+1} + \theta_i )h_x  +\sum_{ \mathcal{V}_h^M(\theta)} M h_x  +\sum_{ {^m\mathcal{V}_h (\theta)}} m h_x - (u^n,1) - \frac{\tau}  {2}  \big( S^{n} + S^{n+1} , 1\big), \\ 
					 & F_n(\xi) :=G_n \left(\tau \, (\xi- \xi^n -\lambda_i^n H^{\prime} (u_i^n ) ) \right).
				\end{aligned}
			\end{align*}
Therefore,  $\xi^{n+1}$ is a solution to the nonlinear algebraic equation $F_n(\xi)=0$. However, since the derivative $F_n^{\prime}(\xi)$ might be infeasible, therefore, 
we employ the simple secant method to find the solution $\xi^{n+1}$:
	\begin{align}\label{xi_secant}
			\xi_{(0)}=0,~ \xi_{(1)}=\malO (\tau), ~	\xi_{(k+1)}=\xi_{(k)}-\frac{F_n(\xi_{(k)})(\xi_{(k)}-\xi_{(k-1)} )}{F_n(\xi_{(k)})-F_n(\xi_{(k-1)})}, ~~ k\ge 1. 
			\end{align}
In our practical simulations, \eqref{xi_secant} often converges in just a few iterations, making the cost negligible  \cite{cheng2022new2}.
			
Once $\xi^{n+1}$ is known, then $\theta_i^{n+1}$ is given by \eqref{lag4}, and thus we can update $(u_i^{n+1}, \lambda_i^{n+1})$ via \eqref{lag5}.	We now summarize the implementation of the BP-MC-HOC-Splitting scheme as Algorithm \ref{alg_1}.
\begin{algorithm} [!ht] \small
			\caption{\small  BP-MC-HOC-Splitting algorithm}\label{alg_1}
			\begin{algorithmic}  [1]
				\STATE \textbf{Input:} Initial values $u_h^0:=u^o$, $\lambda_h^0 :=0$, $\xi^0 := 0$ 			
				\STATE \textbf{Output:} Numerical solutions that satisfy $u_i^{n+1} \in[m, M]$ and $(u_h^{n+1} ,1 )=(u_h^{n} ,1)+ \frac{\tau}  {2}  (S^{n} + S^{n+1}, 1)$.		

               \FOR{$n=0$ to $N_t-1$}
			
				\STATE \textbf{Step 1. (prediction step)}
				 Solve the HOC-Splitting scheme \eqref{scheme:1d:MC:p1}--\eqref{scheme:1d:MC:p4}  to get a prediction solution $\tilde u_h^{n+1}$ with inputs $(u_h^{n}, \lambda_h^{n}, \xi^{n})$.  
				
				\STATE \qquad Compute  $ u_h^{n,1 }  $ by solving the cyclic tridiagonal linear algebraic system \eqref{scheme:1d:MC:p1}  with initial inputs $ u_h^n$ and  $\xi^n$ from  $t_n$ to $t_{n+{1}/{2}}$;
					
				\STATE \qquad Compute $   u_h^{n,2 }  $ by  solving the cyclic tridiagonal linear algebraic system \eqref{scheme:1d:MC:p2}  using  $ u_h^{n,1 } $ and $\xi^n$  from  $t_n$ to $t_{n+1}$;
					
				\STATE \qquad Compute $ u_h^{n,3 } $ by solving the cyclic tridiagonal linear algebraic system \eqref{scheme:1d:MC:p3}  using $ u_h^{n,1}$ and 
    $u_h^{n,2 } $  from  $t_n$ to $t_{n+1}$;
					
				\STATE \qquad Compute $\tilde u_h^{n+1} $ by solving the cyclic tridiagonal linear algebraic system \eqref{scheme:1d:MC:p4}  using  $ u_h^{n,3 } $ and  $\xi^n$ from $t_{n+{1}/{2}}$ to $t_{n+1}$;
					
				\STATE \textbf{Step 2. (correction step)}  
				Correct the obtained solution $\tilde u_h^{n+1} $ to desire $(u_h^{n+1}, \lambda_h^{n+1}, \xi^{n+1})$. 
			
				\STATE \qquad  Compute $\xi^{n+1}$ by the secant method \eqref{xi_secant} with initial guesses $ \xi_{(0)}=0$, $\xi_{(1)}=\malO (\tau)$, inputs $(u_h^{n}, \lambda_h^{n}, \xi^{n})$  and predicted solution $\tilde u_h^{n+1}$;
				\STATE \qquad Update $\theta^{n+1}$ by formula \eqref{lag4};
				\STATE \qquad Judge the values $\{\tilde{u}_i^{n+1}   +\theta_i^{n+1}\}$ to determine $  (u_h^{n+1}, \lambda_h^{n+1} )$   by formula \eqref{lag5}.
		
		 \ENDFOR				
	   \end{algorithmic}
	\end{algorithm}
						
\section{High-order ADI splitting schemes for 2D model}\label{sec:2D:ADI}
In this section, we concern on high-order structure-preserving schemes for the two-dimensional nonlinear convection diffusion equation with periodic boundary conditions: 
\begin{align}  \label{model:2D:e1} 
	\left\{     
	\begin{aligned}
    & u_t   +  f(u)_ x  + g(u)_y  - \gamma^x u_{xx}  -\gamma^{y}  u_{yy}
      = S(\bm x,t),   &\quad 0<  x,y < L,~  0<   t \leq T,   \\
    &  u(\bm x,  0)=u^o(\bm x),  &\quad  0\leq  x,y \leq L, \\
    & u(0, y, t)  =u(L, y, t),   &\quad 0\leq  y \leq L, ~ 0\leq   t \leq T,\\   
    &  u(x,0, t)  =u(x, L, t),  &\quad 0\leq  x \leq L, ~ 0\leq   t \leq T.
  	\end{aligned}
  \right.
\end{align}
 where $\bm x=(x,y)$, $ \gamma^x$ and $\gamma^y$  are two given positive constants, and $f(u)$, $  g(u) $ and $S(\bm x,t)$ are well-defined smooth functions. 
 
Similar to \eqref{strang_1}--\eqref{strang_3}, we can express the second-order operator splitting for model \eqref{model:2D:e1} as finding $u^{n+1}(\bm x) = u^{***}(\bm x,t_{n+1})$ such that
 \begin{numcases}{}
 	u^{*}_{t} -\gamma^{x} u^{*}_{xx}- \gamma^{y} u^{*}_{yy}= 0,   & $t\in [t_n, t_{n+1/2}]; \quad u^{*}(\bm x,{t_n})=u^n(\bm x)$,\label {2d:strang_1} \\
 	u^{**}_{t} + f(u^{**})_x + g(u^{**})_y = S(\bm x,t),  & $t\in [t_n, t_{n+1}]; \quad u^{**}(\bm x,t_n)=u^{*}(\bm x,t_{n+{1}/{2}})$, \label {2d:strang_2} \\
 	u^{***}_{t} -\gamma^x u^{***}_{xx}-\gamma^{y} u^{***}_{yy} = 0,  &  $t\in [t_{n+1/2}, t_{n+1}]; \quad u^{***}(\bm x,t_{n+{1}/{2}})=u^{**}(\bm x,t_{n+1})$, \label {2d:strang_3}
 \end{numcases}	
 for $n = 0,1,\cdots, N_t-1$. 
 
 Then, similar to Sect. \ref{sec:1D}, using the Lagrange multiplier approach and HOC discretization, we can construct similar structure-preserving HOC-Splitting schemes for model \eqref{model:2D:e1} based on operator splitting \eqref{2d:strang_1}--\eqref{2d:strang_3}. However, these directly extended schemes shall still be computationally expensive for the 2D model, since the degree of freedoms increases four times, when the spatial mesh is double refined. This limitation motivates us to develop some efficient high-order accurate ADI schemes for the approximation of \eqref{model:2D:e1} to greatly reduce the computational cost.

\subsection{HOC-ADI-Splitting scheme}

First, we introduce some notations. Similar to the one-dimensional case, we define
 $$\Omega_h=\left\{ \bm x_{i,j}:=(x_i,y_j)  \mid 0 \leq i \leq N_x, ~ 0 \leq j \leq N_y \right\}, ~~h_x = L/N_x, ~~h_y = L/N_y, $$
$$\mathcal{V}_h = \left\{ u = \{u_{i,j}\} \mid u_{i+N_x,j} = u_{i,j}, u_{i,j+N_y} = u_{i,j} \right\},$$
where $N_x$ and  $N_y$ are given positive integers.  
For any grid function $v \in \mathcal{V}_h $, we  define the difference operators  
 $\malD_{\hat{z}} v_{i, j} $, $\delta_z^2  v_{i,j} $, and compact operators $\malA_{z}$ and  $\malB_{z}$, for $z=x$ and $y$, accordingly. For simplicity, denote $\malA:=\malA_{x}\malA_{y}$ and $\malB:=\malB_{x}\malB_{y}$.
 
Next, analogs to the one-dimensional model, a two-dimensional HOC-Splitting scheme for \eqref{model:2D:e1} can be proposed as
 \begin{numcases}  {} 
   \Big(\malA  - \frac{\gamma^x \tau  }{4 } \malA_{y}  \delta_x^2 
      - \frac{\gamma^{y} \tau  }{4 }  \malA_{x}  \delta_y^2\Big)  u_{i,j}^{n,1} 
    =  \Big(   \malA +  \frac{\gamma^x \tau  }{4 } \malA_{y}  \delta_x^2    
          +  \frac{\gamma^{y} \tau  }{4 }  \malA_{x}  \delta_y^2  \Big)u_{i,j}^n,    \label{HOC:2D:e1}\\
    \malB   u_{i,j}^{n,2}  
      =  \malH_1 (  u_{i,j}^{n,1}), ~~
      \malB   u_{i,j}^{n,3}  
    =  \malH_2 ( u_{i,j}^{n,1}, u_{i,j}^{n,2}),  \label{HOC:2D:e3} \\
   \Big(\malA  - \frac{\gamma^x \tau }{4} \malA_{y}  \delta_x^2 - \frac{\gamma^{y} \tau  }{4 }  \malA_{x}  \delta_y^2\Big)  u_{i,j}^{n+1} =   \Big(   \malA +
   \frac{\gamma^x \tau  }{4 } \malA_{y}  \delta_x^2    +  \frac{\gamma^{y} \tau  }{4 }  \malA_{x}  \delta_y^2  \Big)u_{i,j}^{n,3},  \label{HOC:2D:e4}
 \end{numcases}
where $$ \mathcal{H}_1(v ) := \malB v -\tau \malB_{y}   \malD_{\hat{x}}   f(v)-\tau \malB_{x}   \malD_{\hat{y}}   g(v) + \tau \malB  S^n,$$ 	
 $$  \mathcal{H}_2(w,v):=\malB  w  -\frac{\tau}  {2} \malB_{y}   \malD_{\hat{x}}  \left(  f(w)+ f(v)  \right)  -\frac{\tau}{2} \malB_{x} \malD_{\hat{y}}   \left( g(w) + g(v)   \right)  
  +\frac{\tau}  {2} \malB  \big(  S^n + S^{n+1}\big).$$
 
It  is clear that the  dimensions of the resulting coefficient matrices of \eqref{HOC:2D:e1}--\eqref{HOC:2D:e4} are all of order $\malO(N_x N_y)$, so the computational cost has significantly increased, compared to the one-dimensional analogs \eqref{1d:scheme:convec1}--\eqref{1d:scheme:diffu2}. In order to further improve computational efficiency, we shall employ  the D{$^{\prime}$} Yakonov type ADI approach for solving \eqref{HOC:2D:e1}--\eqref{HOC:2D:e4}, which reduces the multidimensional large-scale problem into a set of independent small-scale subproblems of one dimension. To this end, we add a high-order perturbation term
	\begin{align} \label{HOC-ADI:e0}
		\frac{\tau^2}{16} \gamma^x \gamma^y  \delta_x^2  \delta_y^2 (u_{i,j}^{n,1}- u_{i,j}^n)
	\end{align}
 to equation \eqref{HOC:2D:e1} to arrive at
  \begin{align}\label{HOC-ADI:e1}
   \Big( \malA_{x}   -\frac{\gamma^x \tau}{4} {\delta_x^2} \Big)
   \Big( \malA_{y}   -\frac{\gamma^y \tau}{4} {\delta_y^2}  \Big) u_{i,j}^{n,1} 
   = \Big(  \malA_{x}   + \frac{\gamma^x \tau}{4} {\delta_x^2}   \Big)
  \Big(  \malA_{y} + \frac{\gamma^y \tau}{4} {\delta_y^2} \Big)    u_{i,j}^n.	
 \end{align}
Then, the implementation of \eqref{HOC-ADI:e1} reduces the solution of the two-dimensional system to that of two sets of one-dimensional system as below:
 \begin{numcases}{}
	  \Big( \malA_{x}   -\frac{\gamma^x \tau}{4} {\delta_x^2} \Big)   \hat u_{i,j}^{n,1} 	
 =  \Big(  \malA_{x}   + \frac{\gamma^x \tau}{4} {\delta_x^2}   \Big)
	\Big(  \malA_{y} + \frac{\gamma^y \tau}{4} {\delta_y^2} \Big)    u_{i,j}^n, \notag \\
	 \Big( \malA_{y}   -\frac{\gamma^y \tau}{4} {\delta_y^2}  \Big) u_{i,j}^{n,1} 
	 = \hat u_{i,j}^{n,1}.	\notag 
\end{numcases} 
    
Next, note that each equation of \eqref{HOC:2D:e3} can be easily split into two parts as
 	  \begin{equation*}\label{HOC-ADI:e2:0} 
 		 \malB_{x} \hat u_{i,j}^{n,2}  =  \malH_1 (  u_{i,j}^{n,1})   \quad  \textrm{and \quad}
 		   \malB_{y}   u_{i,j}^{n,2}  =\hat  u_{i,j}^{n,2},   
 	\end{equation*}		
 and
 	  \begin{equation*}\label{HOC-ADI:e3:0} 
 		 \malB_{x} \hat u_{i,j}^{n,3}  =   \malH_2  ( u_{i,j}^{n,1},    u_{i,j}^{n,2})   \quad  \textrm{and \quad}
 		 \malB_{y}   u_{i,j}^{n,3}  =\hat u_{i,j}^{n,3}.    
 	\end{equation*}		

 Finally, in the remaining half-step $ [t_{n+{1}/{2}}, t_{n+1}]$, the ADI scheme for \eqref{HOC:2D:e4} can be similarly given
 \begin{numcases}{} 
 	\Big( \malA_{x}   -\frac{\gamma^x \tau}{4} {\delta_x^2} \Big)  \hat u_{i,j}^{n,4} 	
 	= \Big(  \malA_{x}   + \frac{\gamma^x \tau}{4} {\delta_x^2}   \Big)
 	\Big(  \malA_{y} + \frac{\gamma^y \tau}{4} {\delta_y^2} \Big)  u_{i,j}^{n,3}, \notag  \\
 	 \Big( \malA_{y}   -\frac{\gamma^y \tau}{4} {\delta_y^2}  \Big) u_{i,j}^{n+1} = \hat u_{i,j}^{n,4}. \notag
 \end{numcases}		

In summary, the HOC-ADI-Splitting scheme for the 2D model \eqref{model:2D:e1} reads as follows:
  \begin{numcases}{}
	  \Big( \malA_{x}   -\frac{\gamma^x \tau}{4} {\delta_x^2} \Big)   \hat u_{i,j}^{n,1} 	
 =  \Big(  \malA_{x}   + \frac{\gamma^x \tau}{4} {\delta_x^2}   \Big)
	\Big(  \malA_{y} + \frac{\gamma^y \tau}{4} {\delta_y^2} \Big)    u_{i,j}^n, \label{HOC-ADI:e1a}  \\
	 \Big( \malA_{y}   -\frac{\gamma^y \tau}{4} {\delta_y^2}  \Big) u_{i,j}^{n,1} 
	 = \hat u_{i,j}^{n,1},	\label{HOC-ADI:e1b} \\
       \malB_{x} \hat u_{i,j}^{n,2}  =  \malH_1 (  u_{i,j}^{n,1}),   \quad 
       \malB_{y}   u_{i,j}^{n,2}  =\hat  u_{i,j}^{n,2}, \quad 
       	 \malB_{x} \hat u_{i,j}^{n,3}  =   \malH_2  ( u_{i,j}^{n,1},    u_{i,j}^{n,2}),   \quad  
 	 \malB_{y}   u_{i,j}^{n,3}  =\hat u_{i,j}^{n,3},\label{HOC-ADI:e3}  \\
   	\Big( \malA_{x}   -\frac{\gamma^x \tau}{4} {\delta_x^2} \Big)  \hat u_{i,j}^{n,4}	
	= \Big(  \malA_{x}   + \frac{\gamma^x \tau}{4} {\delta_x^2}   \Big)
	\Big(  \malA_{y} + \frac{\gamma^y \tau}{4} {\delta_y^2} \Big)  u_{i,j}^{n,3}, \label{HOC-ADI:e4a} \\
	 \Big( \malA_{y}   -\frac{\gamma^y \tau}{4} {\delta_y^2}  \Big) u_{i,j}^{n+1} = \hat u_{i,j}^{n,4}.	\label{HOC-ADI:e4b}
\end{numcases}

\begin{remark}\label{rem:s4e1} Note that the ADI approach is reduced to solve a series of one-dimensional subproblems, that is, solving the intermediate variables $\hat u_h^{n,1}$ -- $\hat u_h^{n,4}$ along the $x$-direction,  and the intermediate variables  $u_h^{n,1}$ -- $u_h^{n,3}$ and the solution $u_h^{n+1}$ along the $y$-direction. It can be seen that all the runs only require solving some independent one-dimensional constant cyclic tridiagonal systems with dimensions $N_x$ and $N_y$. As pointed out in Remark \ref{rem:e1}, both the memory requirement and computational complexity of \eqref{HOC-ADI:e1a}--\eqref{HOC-ADI:e4b} are only of order $\malO (N_xN_y)$. 
\end{remark}

Define the standard discrete $L_2$ inner product and norm
$$(w,v):= \sum_{i=1}^{N_x} \sum_{j=1}^{N_y} h_x h_y w_{i,j} v_{i,j}, ~~\|v\|:=\sqrt{(v,v)}, \quad v,w \in \mathcal{V}_h.$$    
Similar to Theorem \ref{thm:mass:1d}, we can also prove the discrete mass-conservative property for \eqref{HOC-ADI:e1a}--\eqref{HOC-ADI:e4b}.
\begin{theorem}(Discrete Mass Conservation) \label{thm:mass:2d}
        Let $u_h^n= \{u_{i,j}^{n}\} \in \mathcal{V}_h$ be the solution of the HOC-ADI-Splitting scheme \eqref{HOC-ADI:e1a}--\eqref{HOC-ADI:e4b}. Then there holds 
        \begin{align}\label{MassCons:D:2}
	        ( u_h^{n+1}, 1 ) =( {u}^{o} , 1 )+ \frac{\tau}  {2} \sum_{k=0}^{n}  (S^{k} + S^{k+1}  , 1 ). 
        \end{align}
       \end{theorem}

	 \subsection{BP-HOC-ADI-Splitting scheme}\label{sec_2d_lag_bound}
 As the one-dimensional case, by introducing a bound-preserving Lagrange multiplier $\lambda ^{n+1}_{i,j}$, and employing the developed HOC-ADI-Splitting scheme \eqref{HOC-ADI:e1a}--\eqref{HOC-ADI:e4b}, we can construct the following BP-HOC-ADI-Splitting scheme, which preserves the bounds of the physical quantity but may not necessarily guarantee the mass conservation. 
  \paragraph{Step 1 (Predictor).} 	 Solve $\overline  {u}_h^{n+1}$ from $u_h^{n}$ using the  proposed scheme \eqref{HOC-ADI:e1a}--\eqref{HOC-ADI:e4b} with  modifications in  \eqref{HOC-ADI:e1a} and  \eqref{HOC-ADI:e4a}  that	
	\begin{numcases} {}  
		\Big( \malA_{x}   -\frac{\gamma^x \tau}{4} {\delta_x^2} \Big) \hat u_{i,j}^{n,1} 	  =  \Big(  \malA_{x}   + \frac{\gamma^x \tau}{4} {\delta_x^2}   \Big)
				\Big(  \malA_{y} + \frac{\gamma^y \tau}{4} {\delta_y^2} \Big)  u_{i,j}^n + \frac{\tau}{2}\malA (\lambda_{i,j}^n  H^{\prime}(u_{i,j}^n)),  \label{HOC-ADI:BP:e1a}\\
				 \Big( \malA_{y}   -\frac{\gamma^y \tau}{4} {\delta_y^2}  \Big) u_{i,j}^{n,1} = \hat u_{i,j}^{n,1},\label{HOC-ADI:BP:e1b}	\\
		       \malB_{x} \hat u_{i,j}^{n,2}  =  \malH_1 (u_{i,j}^{n,1}), \quad
               \malB_{y}   u_{i,j}^{n,2}  =\hat  u_{i,j}^{n,2}, \quad				
			 \malB_{x}   \hat {u}_{i,j}^{n,3}  
				=    \malH_2  (u_{i,j}^{n,1},    u_{i,j}^{n,2}),	\quad 		
              \malB_{y} u_{i,j}^{n,3}=  \hat {u}_{i,j}^{n,3}, \label{HOC-ADI:BP:e1d} \\
			 	\Big( \malA_{x}   -\frac{\gamma^x \tau}{4} {\delta_x^2} \Big)  \hat u_{i,j}^{n,4} 	
				= \Big(  \malA_{x}   + \frac{\gamma^x \tau}{4} {\delta_x^2}   \Big)
				\Big(  \malA_{y} + \frac{\gamma^y \tau}{4} {\delta_y^2} \Big)  u_{i,j}^{n,3} + \frac{\tau}{2}\malA (\lambda_{i,j}^n  H^{\prime}(u_{i,j}^n) ), \label{HOC-ADI:BP:e1e}\\
			 \Big( \malA_{y}   -\frac{\gamma^y \tau}{4} {\delta_y^2}  \Big) \overline u_{i,j}^{n+1} = \hat u_{i,j}^{n,4}.\label{HOC-ADI:BP:e1f}
		\end{numcases}
		
		\paragraph{Step 2 (Corrector).}  Solve $(u_h^{n+1}, \lambda_h^{n+1})$ via the cut-off approach
		\begin{align}\label{HOC-ADI:BP:e2} 
			\big(u_{i,j}^{n+1} , \lambda_{i,j}^{n+1} \big)
			 =\left\{\begin{array}{cll}
				\big(\overline{u}_{i,j}^{n+1} +\eta_{i,j}^{n+1} , 0 \big), & \text { if } \quad m< \overline{u}_{i,j}^{n+1} +\eta_{i,j}^{n+1} <M,   \\[0.1in]
				\Big( m, \frac{ m-( \overline{u}_{i,j}^{n+1} +\eta_{i,j}^{n+1} )}{\tau H^{\prime}(m)}\Big), & \text { if } \quad \overline{u}_{i,j}^{n+1} +\eta_{i,j}^{n+1}  \leq m, \\[0.1in]
				\Big(M, \frac{  M-( \overline{u}_{i,j}^{n+1} +\eta_{i,j}^{n+1} ) }{ \tau H^{\prime}(M)}\Big), & \text { if} \quad \overline{u}_{i,j}^{n+1} +\eta_{i,j}^{n+1}  \geq M,
			\end{array}  \right. 
		\end{align}
			where 
		$ \eta_{i,j}^{n+1} :=- { \tau} \lambda_{i,j}^n H^{\prime}(u_{i,j}^n)$.

\subsection{ BP-MC-HOC-ADI-Splitting scheme}\label{sec_2d_lag_mass}
In order to further ensure mass conservation, we introduce another Lagrange multiplier $\xi^{n+1}$, which is independent of spatial variables, to enforce the mass conservation in the correction step. The new scheme, named BP-MC-HOC-ADI-Splitting, is proposed as follows:	
	\paragraph{Step 1 (Predictor).}  Solve $\tilde {u}_h^{n+1}$ from $u_h^{n}$  using the given scheme 
		 \begin{numcases}  {} 
			 \Big( \malA_{x}   -\frac{\gamma^x \tau}{4} {\delta_x^2}\Big)   \hat u_{i,j}^{n,1} 	
				=  \Big(  \malA_{x}   + \frac{\gamma^x \tau}{4} {\delta_x^2}  \Big)
			    \Big(  \malA_{y} + \frac{\gamma^y \tau}{4} {\delta_y^2}\Big)  u_{i,j}^n + \frac{\tau}{2} 
			    \malA   \Big(  \xi^{n} + \lambda_{i,j}^n  H^{\prime}(u_{i,j}^n )\Big),  \label{HOC-ADI:Mass:e1} \\
			  \Big( \malA_{y}   -\frac{\gamma^y \tau}{4} {\delta_y^2} \Big) u_{i,j}^{n,1} = \hat u_{i,j}^{n,1}, \label{HOC-ADI:Mass:e2}	\\
			  \malB_{x} \hat u_{i,j}^{n,2}  =  \malH_1 (  u_{i,j}^{n,1}), \quad 
                \malB_{y}   u_{i,j}^{n,2}  =\hat  u_{i,j}^{n,2}, \quad			
			  \malB_{x}   \hat {u}_{i,j}^{n,3}  
				=   \malH_2  ( u_{i,j}^{n,1},    u_{i,j}^{n,2}), \quad 	 
               \malB_{y} u_{i,j}^{n,3}=  \hat {u}_{i,j}^{n,3},  \label{HOC-ADI:Mass:e4}\\
			  \Big( \malA_{x}   -\frac{\gamma^x \tau}{4} {\delta_x^2}\Big)  
     \hat u_{i,j}^{n,4} 	
				 = \Big(  \malA_{x}   + \frac{\gamma^x \tau}{4} {\delta_x^2}  \Big)
				\Big(  \malA_{y} + \frac{\gamma^y \tau}{4} {\delta_y^2}\Big)  u_{i,j}^{n,3} + \frac{\tau}{2} 
				\malA   \Big(  \xi^{n} + \lambda_{i,j}^n  H^{\prime}(u_{i,j}^n )\Big), \label{HOC-ADI:Mass:e5}\\
			  \Big( \malA_{y}   -\frac{\gamma^y \tau}{4} {\delta_y^2} \Big)  \tilde u_{i,j}^{n+1} = \hat u_{i,j}^{n,4}.\label{HOC-ADI:Mass:e6}
		\end{numcases}
			
\paragraph{Step 2 (Corrector).}   Solve $(u_h^{n+1}, \lambda_h^{n+1}, \xi^{n+1})$ via the cut-off approach
			\begin{align}\label{HOC-ADI:Mass:correct}  
				& (u_{i,j}^{n+1}   , \lambda_{i,j}^{n+1} )
                 =\left\{\begin{array}{ll}
					\big(\tilde{u}_{i,j}^{n+1}   +\theta_{i,j}^{n+1}   , 0\big), & \text { if } \quad m<\tilde{u}_{i,j}^{n+1}   +\theta_{i,j}^{n+1}   <M, \\
					\Big(m, \frac{  m-(\tilde{u}_{i,j}^{n+1}   +\theta_{i,j}^{n+1}   )}{\tau H^{\prime}(m)}\Big), & \text{ if } \quad \tilde{u}_{i,j}^{n+1}   +\theta_{i,j}^{n+1}  \leq m, \\
					\Big(M, \frac{ M-(\tilde{u}_{i,j}^{n+1}   +\theta_{i,j}^{n+1} )}{{\tau } H^{\prime}(M)}\Big), & \text { if } \quad \tilde{u}_{i,j}^{n+1}   +\theta_{i,j}^{n+1}  \geq M,
				\end{array}  \right. 
			\end{align}
where $\theta_{i,j}^{n+1} := {\tau}\big( \xi^{n+1}- \xi^n - \lambda_{i,j}^n H^{\prime}(u_{i,j}^n)\big) $, and similarly, the multiplier $\xi^{n+1}$ is solved iteratively via the secant method \eqref{xi_secant}.
				
\begin{remark} \label{rem:2d:matrix}
The coefficient matrices of the proposed scheme \eqref{HOC-ADI:Mass:e1}--\eqref{HOC-ADI:Mass:e6} are totally the same as the one-dimensional case and only need to be generated once.  Therefore, the implementation of the ADI approach is basically the same as the one-dimensional BP-MC-HOC-Splitting scheme along two separate spatial directions.			
\end{remark}	

 \begin{remark} 
It is clear that all the 
 schemes constructed in subsections \ref{subsec:1d:lag}, \ref{sec_2d_lag_bound} and \ref{sec_2d_lag_mass} can automatically ensure the $L_\infty$  bound for $ u_h^n$, that is, the proposed structure-preserving HOC-Splitting schemes are naturally stable. However, the error analysis for the BP-MC-HOC-ADI-Splitting scheme is still challenging, and we only present error estimates for the BP-HOC-ADI-Splitting schemes with cut-off approach \cite{cheng2022new2}, see Sect. \ref{sec:ErrEst} for details. The one-dimensional case is simpler compared to the two-dimensional case, and the interested readers are referred to prove it in a similar way.
 \end{remark}
 
 %

We now summarize the above BP-MC-HOC-ADI-Splitting scheme \eqref{HOC-ADI:Mass:e1}--\eqref{HOC-ADI:Mass:correct} into Algorithm \ref{alg_2}.
\begin{algorithm}  [!ht] \small
	\caption{\small BP-MC-HOC-ADI-Splitting algorithm }
	\label{alg_2}
	\begin{algorithmic}  [1]
		\STATE \textbf{Input:} Initial values $u_h^{0} := u^o$,  $\lambda_h^0 :=0$, $\xi^0 := 0$ 
		
		\STATE \textbf{Output:} Numerical solutions that satisfy $ u^n_{i,j}\in[m, M]$ and $(u_h^{n+1} ,1 )=(  u_h^{n} ,1)+ \frac{\tau}  {2}  (S^{n} + S^{n+1}, 1)$.
		
		\FOR{$n=0$ to $N_t-1$}
		
		\STATE \textbf{Step 1.} (\textbf{prediction step})
		Solve the HOC-ADI-Splitting scheme \eqref{HOC-ADI:Mass:e1}--\eqref{HOC-ADI:Mass:e6}  to  get a prediction solution $\tilde {u}_h^{n+1}$ with inputs $(u_h^{n}, \lambda_h^{n}, \xi^{n})$ .
		
		\STATE \qquad Compute  $\hat u_h^{n,1} $ along $x$-direction by \eqref{HOC-ADI:Mass:e1} with inputs $(u_h^{n}, \lambda_h^{n}, \xi^{n})$; 
		
		\STATE  \qquad Compute $   u_h^{n,1}  $ along $y$-direction  by \eqref{HOC-ADI:Mass:e2} using  $ \hat u_h^{n,1} $; 
		
		\STATE \qquad Compute $  \hat u_h^{n,2 } $ along $x$-direction  using $ u_h^{n,1} $, and then compute $ u_h^{n,2 } $ along $y$-direction using  $\hat u_h^{n,2 } $  by the first two equations of \eqref{HOC-ADI:Mass:e4};
			
		\STATE \qquad Compute $  \hat u_h^{n,3 } $ along $x$-direction using  $ \{u_h^{n,1 }, u_h^{n,2 }\}$, and then compute $  u_h^{n,3} $ along $y$-direction using  $\hat u_h^{n,3 } $   by  the last two equations of  \eqref{HOC-ADI:Mass:e4};
		
		\STATE \qquad Compute $\hat u_h^{n,4} $ along $x$-direction by \eqref{HOC-ADI:Mass:e5} using $ \{ u_h^{n,3 }, \xi^n,   \lambda_h^n \}$;  
		
		\STATE \qquad  Compute $\tilde u_h^{n+1} $ along $y$-direction by \eqref{HOC-ADI:Mass:e6} using  $\hat u_h^{n,4} $. 
		
		\STATE \textbf{Step 2.} (\textbf{correction step})  {Correct the   estimated value $\tilde u_h^{n+1}$ to desire the solution $u_h^{n+1}$.}
		
		\STATE \qquad  Compute $\xi^{n+1}$ by the secant method \eqref {xi_secant} with initial guesses $ \xi_{(0)}=0$, $\xi_{(1)}=\malO (\tau)$, inputs $(u_h^{n}, \lambda_h^{n}, \xi^{n})$  and predicted solution $\tilde u_h^{n+1}$;
		
		\STATE \qquad Update $\theta^{n+1}$ via the formula $\theta_{i,j}^{n+1} = \tau\,  (\xi^{n+1}- \xi^n - \lambda_{i,j}^n H^{\prime}(u_{i,j}^n))$;
		
		\STATE \qquad Judge the values $\{\tilde{u}_{i,j}^{n+1}   +\theta_{i,j}^{n+1}\}$ to determine $(u_h^{n+1}, \lambda_h^{n+1} )$   by  formula \eqref{HOC-ADI:Mass:correct}.
		
		\ENDFOR
	\end{algorithmic}
\end{algorithm}

\section{Error analysis of the bound-preserving scheme}\label{sec:ErrEst}
 In this section, we aim to present an error analysis for   the BP-HOC-ADI-Splitting scheme \eqref{HOC-ADI:BP:e1a}--\eqref{HOC-ADI:BP:e2} in discrete $L_2$ norm by analyzing each subprocess separately. Without loss of generality, we assume $S=0$. 
 As mentioned in Remark \ref{rm_lambda_0}, the Lagrange multiplier $\lambda$ provides an alternative interpretation of the high-order cut-off approach, and in actual computation we can replace $\lambda_{i,j}^n  H^{\prime}(u_{i,j}^n )$ by zero in \eqref{HOC-ADI:BP:e1a}--\eqref{HOC-ADI:BP:e2}, which will not affect the accuracy \cite{cheng2022new2}. 
Thus, in the following, we only pay attention to error analysis of the bound-preserving scheme \eqref{HOC-ADI:BP:e1a}--\eqref{HOC-ADI:BP:e2} with $\lambda_{i,j}^n  H^{\prime}(u_{i,j}^n )$ replaced by zero, which is easier to analyze. Throughout this section, we use $C$ to denote a generic positive constant that may depend on the given data but is independent of the mesh parameters.

To begin with, we introduce some useful lemmas.
 \begin{lemma}[\cite{liao2010maximum}] \label{lem:err:Taylor_HOC} 
 Let $v(x) \in {C}^6[0,1]$ and $\zeta(\lambda)=5(1-\lambda)^3-3(1-\lambda)^5$, then
 	$$
 	\begin{aligned}
 		& \frac{v^{\prime \prime}\left(x_{i+1}\right)+10 v^{\prime \prime}\left(x_i\right)+v^{\prime \prime}\left(x_{i-1}\right)}{12}=\frac{v\left(x_{i+1}\right)-2 v\left(x_i\right)+v\left(x_{i-1}\right)}{h_x^2} \\
 		& \quad +\frac{h_x^4}{360} \int_0^1\left[v^{(6)}\left(x_i-\lambda h_x\right)+v^{(6)}\left(x_i+\lambda h_x\right)\right] \zeta(\lambda) d \lambda, \quad 1 \leq i \leq N_x.
 	\end{aligned}
 	$$	
 \end{lemma}
 \begin{lemma}[\cite{xu2024maximum}] \label{lem:err:AB_bound}
 For any grid function $v \in \mathcal{V}_h$ and $z=x,y$, it holds that  
  \begin{align*}
    & \frac{2}{3}\|v\| \le \|\malA_z v\| \leq\|v\|, \quad \frac{2}{3}\|v\|^2 \le (\malA_z v, v) \leq\|v\|^2,
    \quad \frac{4}{9}\|v\| \leq \|\malA v\| \leq\|v\|, 
    \\
    & \frac{1}{3}\|v\| \le \|\malB_z v\| \leq\|v\|, \quad \frac{1}{3}\|v\|^2 \le (\malB_z v, v) \leq\|v\|^2,
    \quad \frac{1}{9}\|v\| \leq \|\malB v\| \leq\|v\|. 
 \end{align*}
 
 \end{lemma}

Lemma \ref{lem:err:AB_bound} shows that the operator $\malA $ is symmetric and positive definite. Thus, we can define the following discrete inner product and associated norm
$$
(v, w)_\malA :=(\malA  v, w), \quad \|v\|_\malA :=\sqrt{(v, v)_\malA}.
$$
It is clear that the two norms $\|\cdot\|_\malA $ and $\|\cdot\|$ are equivalent. Similarly, $ \|\cdot\|_{\malA_x}$ and $ \|\cdot\|_{\malA_y} $ can also be given and all of them are equivalent.
  
Next, we present the main conclusion for the bound-preserving scheme in the following theorem. For simplicity, we write $u (t_{n})$, $u^{n}$ and ${u}^{n}_h$ respectively denote the exact solution, the operator splitting solution and the numerical solution at $\{\bm x_{i,j}\}$. 
\begin{theorem}\label{th:err:BP-HOC-ADI}
Let $ u(\bm x, t) \in  {C}^{6,4}({\Omega} \times[0, T])$ be the solution of model \eqref{model:2D:e1} and $\{u_h^{n} \mid  0 \leq n \leq N_t\big\}$ be the numerical solution of the BP-HOC-ADI-Splitting scheme \eqref{HOC-ADI:BP:e1a}--\eqref{HOC-ADI:BP:e2}.  Suppose that  $f(u)$ and $g(u)$ are sufficiently smooth and  $f^\prime(u), g^\prime(u)$  are Lipschitz continuous with Lipschitz constant $L$. Then, for any sufficiently small
$\tau$, $h_x$ and $h_y$, 
  satisfying the CFL condition  
\begin{align} \label{CFL:e1}   
 \max{|f^{\prime}(u)|}\frac{\tau}{h_x} +  \max{|g^{\prime}(u)|}\frac{\tau}{h_y} \le C_0,  
\end{align}   
 it holds
$$
\| u(t_{n})- {u}^{n}_h \| \leq C\big(\tau^2+ h_x^4+h_y^4\big), \quad 1 \leq n \leq N_t,
$$
where $C_0$ is the positive CFL constant independent of the mesh parameters.
\end{theorem} 

The proof of this theorem is divided into several parts. Firstly, the following lemma states the second-order splitting error of the splitting scheme \eqref{2d:strang_1}--\eqref{2d:strang_3}.
 \begin{lemma} [\cite{bookADR}, Chapter IV]\label{lem:err:strang}
 	Let $u(\bm x,  t)$ be the solution of  model \eqref{model:2D:e1} and $\left\{u^n(\bm x) \mid 0 \leq n \leq N_t\right\}$ be the solution of the operator splitting scheme \eqref{2d:strang_1}--\eqref{2d:strang_3}. Then there holds
 	$$
 	\|u(t_n)-u^n\| \leq C \tau^2, \quad 1 \leq n \leq N_t.
 	$$
 \end{lemma} 
 
Next, we present the error estimate for $e^{n}:=u^n - {u}^{n}_h$ in three steps according to the three processes of the operator splitting approach. Firstly, we estimate the diffusion process \eqref{HOC-ADI:BP:e1a}--\eqref{HOC-ADI:BP:e1b} in the first half time step $[t_n,t_{n+1/2}]$. Denote  $e^{*,n} :=  u^*(t_{n+1/2})-u^{n,1}_h$. 
\begin{lemma} \label{lem:err:BP-HOC-ADI-Split:diffusion:1}
   Let $ u^*(\bm x, t)$ be the solution of the diffusion process \eqref{2d:strang_1} and $ \{ u_h^{n,1}  \mid  0 \leq n \leq N_t  \}$ be the numerical solution of the scheme \eqref{HOC-ADI:BP:e1a}--\eqref{HOC-ADI:BP:e1b}. Then, under the assumptions of Theorem \ref{th:err:BP-HOC-ADI}, there holds
$$
\| u^*(t_{n+1/2})- u_h^{n,1} \|_\malA
\leq
 \|u^n- u_h^{n}\|_\malA +C \tau\big(\tau^2+ h_x^4+h_y^4\big), \quad  0 \leq n \leq N_t-1.
$$
\end{lemma}
\begin{proof}
    For convenience, denote $\delta_t u_h^{n,\frac{1}{2}}=(u_h^{n,1}-u_h^n)/{\tau}$ and  $u_h^{n,\frac{1}{2}}=(u_h^{n,1}+u_h^n)/2$. Then, we rewrite \eqref{HOC-ADI:BP:e1a}--\eqref{HOC-ADI:BP:e1b} as follows:
    \begin{align} \label{err:BP:eq_ustar}
        \malA \delta_t u^{n,\frac{1}{2}}_{i,j} - \frac{\gamma^x  }{2}\malA_y \delta_x^2 u^{n,\frac{1}{2}}_{i,j} -\frac{\gamma^y  }{2}\malA_x \delta_y^2 u^{n,\frac{1}{2}}_{i,j} +
        \frac{\gamma^x \gamma^y \tau }{8}\malA \delta_x^2 \delta_y^2 u^{n,\frac{1}{2}}_{i,j} =0.
    \end{align} 
   Besides, define $\delta_t u^{*,\frac{1}{2}}=(u^*(t_{n+1/2})-u^n)/{\tau}$ and  $u^{*,\frac{1}{2}} =(u^*(t_{n+1/2}) +u^n)/2$. It follows from Taylor expansion and Lemma \ref{lem:err:Taylor_HOC} that $ u^{*,\frac{1}{2}}$ satisfies 
    \begin{equation} \label{err:BP:eq_uhat}
        \malA \delta_t u^{*,\frac{1}{2}}(\bm x_{i,j})- \frac{\gamma^x  }{2}\malA_y \delta_x^2 u^{*,\frac{1}{2}}(\bm x_{i,j})  -\frac{\gamma^y  }{2}\malA_x \delta_y^2  u^{*,\frac{1}{2}}(\bm x_{i,j})  +
        \frac{\gamma^x \gamma^y \tau }{8}\malA \delta_x^2 \delta_y^2  u^{*,\frac{1}{2}}(\bm x_{i,j})  =\tau R^{*,\frac{1}{2}}_{i,j},
    \end{equation}
    where $ R^{*,\frac{1}{2}}_{i,j} = \malO \big(\tau^2 +   h_x^4+   h_y^4 \big)$.
   
  Let $e^{*,\frac{1}{2}}: = u^{*,\frac{1}{2}}-u^{n,\frac{1}{2}}_{h}$. Then, subtracting \eqref{err:BP:eq_uhat} from \eqref{err:BP:eq_ustar}, we get the error equation as follows:
    \begin{align} \label{err:BP:eq_ehat}
        \malA \delta_t e^{*,\frac{1}{2}}_{i,j}- \frac{\gamma^x  }{2}\malA_y \delta_x^2 e^{*,\frac{1}{2}}_{i,j}  -\frac{\gamma^y  }{2}\malA_x \delta_y^2  e^{*,\frac{1}{2}}_{i,j}  +
        \frac{\gamma^x \gamma^y \tau }{8}\malA \delta_x^2 \delta_y^2  e^{*,\frac{1}{2}}_{i,j}  = \tau R^{*,\frac{1}{2}}_{i,j}.
    \end{align}

Taking the inner product of \eqref{err:BP:eq_ehat} with $e^{*,\frac{1}{2}}$, combining Lemma \ref{lem:err:AB_bound} and summation by parts, we obtain
    \begin{equation} \label{err:BP:eq_ehat_L2}
        \frac{ \|e^{*,n}\|_\malA^2 - \|e^n\|_\malA^2}{2\tau}+ \frac{\gamma^x  }{2}\|\delta_x e^{*,\frac{1}{2}} \|^2_{\malA_y} +\frac{\gamma^y  }{2} \|\delta_y   e^{*,\frac{1}{2}}\|^2_{\malA_x} +
        \frac{\gamma^x \gamma^y \tau }{8}\|\delta_x \delta_y e^{*,\frac{1}{2}}\|^2_\malA  \le \tau (R^{*,\frac{1}{2}},e^{*,\frac{1}{2}}),
    \end{equation}
which directly leads to
    \begin{align*} 
       \frac{1}{2\tau}\big(\|e^{*,n}\|_\malA^2 - \|e^n\|_\malA^2\big)
        \le \sqrt{\frac{3}{2}} \|R^{*,\frac{1}{2}}\| \frac{\|e^{*,n}\|_\malA+\|e^n\|_\malA}{2} \Longrightarrow 
       \|e^{*,n}\|_\malA \le \|e^n\|_\malA  
        + \sqrt{\frac{3}{2}} \tau \|R^{*,\frac{1}{2}}\|,
    \end{align*}
which completes the proof. \qed
\end{proof}

Secondly, using similar arguments as in Lemma \ref{lem:err:BP-HOC-ADI-Split:diffusion:1}, we can estimate the diffusion process \eqref{HOC-ADI:BP:e1e}--\eqref{HOC-ADI:BP:e1f} in the second half time step $[t_{n+1/2},t_{n+1}]$. Denote $\overline e^{n+1}:=  u^{n+1}-\overline u_h^{n+1}$.

 \begin{lemma} \label{lem:err:BP-HOC-ADI-Split:diffusion:2}    
 Let $ u^{***}(\bm x, t)$ be the solution of the diffusion process \eqref{2d:strang_3} and $ \{\overline u_h^{n+1}  \mid  0 \leq n \leq N_t-1  \}$ be the numerical solution of the scheme \eqref{HOC-ADI:BP:e1e}--\eqref{HOC-ADI:BP:e1f}. Then,  under the assumptions of Theorem \ref{th:err:BP-HOC-ADI}, there holds
$$
\| u^{n+1}- \overline u_h^{n+1} \|_\malA \leq
   \big\|u^{**}(t_{n+1})- u_h^{n,3} \big\|_\malA +C \tau\big(\tau^2+ h_x^4+h_y^4\big), \quad  0 \leq n \leq N_t-1.
$$    
\end{lemma}

Finally, we consider the error estimate of the second subprocess: the transport equation \eqref{HOC-ADI:BP:e1d}  in the entire time step $[t_n,t_{n+1}]$. Denote $e^{**,n} :=  u^{**}(t_{n+1})-u^{n,3}_h$. 
\begin{lemma}  \label{lem:err:BP-HOC-ADI-Split:hyperbolic}
 Let $ u^{**}(\bm x, t)$ be the solution of the convection process \eqref{2d:strang_2} and  $\big\{ {u_h^{n,3}} \mid 0 \leq n \leq N_t\big\}$ be the numerical solution of the scheme \eqref{HOC-ADI:BP:e1d}. Then, under the assumptions of Theorem \ref{th:err:BP-HOC-ADI}, there holds
$$
 \| u^{**}(t_{n+1})-u_h^{n,3} \|_\malA 
 \leq \big(1+C\tau\big) \|  u^*(t_{n+1/2}-u_h^{n,1} \|_\malA+C \tau \big(\tau^2 +  h_x^4+  h_y^4 \big), \quad 0 \leq n \leq N_t-1.
$$
\end{lemma} 
\begin{proof}
Firstly, notice from \eqref{Split:semi:e1} that the SSP-RK2 time-stepping scheme \eqref{HOC-ADI:BP:e1d} can be expressed as
\begin{align}  \label{err:hyper:num}
	\begin{aligned}
		   u^{n,2}_{i,j} =  u^{n,1}_{i,j}  -\tau \malF(u^{n,1}_{i,j}),\quad
		   u^{n,3}_{i,j} = \frac{1}{2}  u^{n,1}_{i,j} + \frac{ 1}{2}  u^{n,2}_{i,j} - \frac{\tau}{2} \malF (u^{n,2}_{i,j}),
	\end{aligned}
\end{align}
where 
$$  \malF\left(v\right) := -\malB_{x}^{-1}  \malD_{\hat{x}} f(v)- \malB_{y}^{-1}    \malD_{\hat{y}} g(v). 
$$ 
Correspondingly, the exact solution $u^{**}(\bm x_{i,j}, t_{n+1})$  satisfies
\begin{align} \label{err:hyper:ext}
	\left\{\begin{aligned}
		&  \check{u}(\bm x_{i,j}, t_{n+1}) :=u^*(\bm x_{i,j}, t_{n+1/2})- \tau \malF\left(u^*(\bm x_{i,j}, t_{n+1/2} )\right), \\
		& u^{**}(\bm x_{i,j}, t_{n+1}) =\frac{1}{2} u^*(\bm x_{i,j}, t_{n+1/2})
            +\frac{1}{2} \check{u}(\bm x_{i,j}, t_{n+1}) -  \frac{\tau}{2}  \malF (\check{u}(\bm x_{i,j}, t_{n+1}))  + \tau R_{i,j}^{**,n},
	\end{aligned}\right.  
\end{align}
 where $ R^{**,n}_{i,j} = \malO \big(\tau^2 +   h_x^4+   h_y^4 \big)$.
 
Denote $\check e^{n}:= \check{u}(t_{n+1}) -u^{n,2}_h$. Then, subtracting \eqref{err:hyper:ext} from \eqref{err:hyper:num}, we get the following error equation
\begin{numcases}{}
   \check e_{i,j}^{n} = e^{*,n}_{i,j} -\tau\big(\malF(u^*(\bm x_{i,j}, t_{n+1/2} )) - \malF(u^{n,1}_{i,j} )\big), \label{err:hyper:err:1}\\
   e^{**,n}_{i,j} =\frac{1}{2} e^{*,n}_{i,j} +\frac{1}{2}  \check e_{i,j}^{n} - \frac{\tau}{2} \big(\malF ( \check{u}^{n}_{i,j} ) -\malF(u^{n,2}_{i,j} )\big) +  \tau R^{**,n}_{i,j}. \label{err:hyper:err:2}
\end{numcases}

Taking the inner product of  \eqref{err:hyper:err:1} with $\malA \check e^n$ and \eqref{err:hyper:err:2} with $\malA  e^{**,n}$, respectively, and then by Cauchy-Schwarz inequality we obtain
\begin{equation}  \label{err:hyper:e_star2} 
  \begin{aligned}
    & \| \check e^n  \|_\malA  \le \|e^{*, n}\|_\malA   + \tau\|\malF(u^*( t_{n+1/2} )) - \malF(u^{n,1}_h )\|_\malA,      \\
    &\|e^{**, n}\|_\malA \le \frac{1}{2} \|e^{*, n}\|_\malA +\frac{1}{2} \|\check e^n\|_\malA + \frac{\tau}{2}\|\malF ( \check{u}^{n}_h ) -\malF(u^{n,2}_h)\|_\malA   +  \tau \|R^{**}\|_\malA.  
 \end{aligned}
\end{equation}

Note that Lemma \ref{lem:err:AB_bound} shows that $\malB_{z}^{-1}$ is bounded, and the Lipschitz continuous condition of $f^{\prime}$ and $g^{\prime}$ implies that
\begin{equation}\label{err:F}
\begin{aligned} 
	|\malF(v_1)- \malF(v_2)|
	&\le |f^{\prime}(v_1)- f^{\prime}(v_2)+ \malO (h_x^4)|+ |g^{\prime}(v_1)- g^{\prime}(v_2)+ \malO (h_y^4)|\\
      & \le L|v_1-v_2|+\malO (h_x^4+h_y^4), \quad v_1,v_2 \in \mathcal{V}_h.
\end{aligned} 
\end{equation}
Therefore, applying \eqref{err:F} into \eqref{err:hyper:e_star2},  under the CFL condition \eqref{CFL:e1}, we have
\begin{equation} \label{err:hyper:e_check} 
  \begin{aligned}
      \| \check e^n  \|_\malA  & \le \|e^{*, n}\|_\malA   + \tau L\|e^{*, n}\|_\malA 
      = (1+\tau L) \|e^{*, n}\|_\malA , 
\end{aligned} 
\end{equation}
and 
   \begin{align*}    
      \|e^{**, n}\|_\malA  \le \frac{1}{2} \|e^{*, n}\|_\malA+ \frac{1+\tau L}{2} \|\check e^n\|_\malA +  \tau \|R^{**}\|_\malA
        \le (1+C\tau) \|e^{*, n}\|_\malA +C \tau \big(\tau^2 +  h_x^4+  h_y^4 \big),
 \end{align*}
 which proves the conclusion. \qed
\end{proof}
 
We are now in a position to prove Theorem \ref{th:err:BP-HOC-ADI}. 
\begin{proof}
 Lemma \ref{lem:err:strang} indicates that, in order to obtain the final error estimate, it suffices to show that 
\begin{align} \label{err:BP:uhat_un}
    \| u^n-u^n_h\| \leq C \big(\tau^2+ h_x^4+h_y^4\big), \quad 1 \leq n \leq N_t.
\end{align}
 
First, it follows from Lemmas \ref{lem:err:BP-HOC-ADI-Split:diffusion:1}--\ref{lem:err:BP-HOC-ADI-Split:hyperbolic} that
\begin{align}  \label{err:bare_n+1:L2}
\big\|\overline e^{n+1}\big\|_\malA =\| u^{n+1}- \overline u_h^{n+1} \|_\malA \leq {(1+ C\tau )}\left\|e^n\right\|_\malA+ {\malO \big(\tau^3 + \tau h_x^4+ \tau h_y^4 \big)}, \quad 0 \leq n \leq N_t-1.
\end{align}

Next, we provide an error estimate for the correction step. For the convenience of analysis, like the 1D formula \eqref{scheme:1d:BP:cor}, we rewrite the correction step \eqref{HOC-ADI:BP:e2}   as follows:     
\begin{align}\label{err:BP-HOC-ADI:cor}
	\left\{   
	\begin{aligned}
		&	\frac{ u_{i,j}^{n+1} -\overline{u}_{i,j}^{n+1}  }{\tau } =  \lambda_{i,j}^{n+1} H^{\prime}(u_{i,j}^{n+1} ), \\
		&	\lambda_{i,j}^{n+1}  \geq 0, ~H(u_{i,j}^{n+1}) \geq 0, ~ \lambda_{i,j}^{n+1}  H(u_{i,j}^{n+1})=0.
	\end{aligned}
	\right.
\end{align}
Then, we derive from \eqref{err:BP-HOC-ADI:cor} that
\begin{align} \label{err:ehat}
   e^{n+1} - \tau\rho_h^{n+1} = \overline e^{n+1}, \quad \rho_h^{n+1}:=-\lambda_h^{n+1} H^\prime(u_h^{n+1}).
\end{align}
Taking the discrete inner product of \eqref{err:ehat} with   $\malA e^{n+1}$, we have
\begin{align}\label{err:l2}
    \| e^{n+1}\|^2_\malA   -\tau(e^{n+1}, \rho_h^{n+1})_\malA = (\overline e^{n+1}, \malA e^{n+1}).
\end{align}

Thanks to the KKT condition of \eqref{err:BP-HOC-ADI:cor}, we obtain
\begin{equation} \label{err:multiplier}
   \begin{aligned}
    -\tau(e^{n+1}, \rho_h^{n+1})_\malA
     & = \tau \big(u^{n+1}- u_h^{n+1}, \lambda_h^{n+1}  H^{\prime}(u_h^{n+1})\big)_\malA + \tau \big(\lambda_h^{n+1},H(  u_h^{n+1})\big)_\malA \\
    & =\tau \big(\lambda_h^{n+1},(e^{n+1})^2\big)_\malA + \tau\big(\lambda_h^{n+1},H(u^{n+1})\big)_\malA \ge 0,
    \end{aligned}
\end{equation}
as  $\lambda_h^{n+1}  \geq 0$ and  $m \leq u^{n+1} \leq M$.
Thus, substituting \eqref{err:multiplier} into \eqref{err:l2}, and using \eqref{err:bare_n+1:L2} we obtain
\begin{align} \label{err:enplus1_end} 
    \| e^{n+1}\|_\malA  \leq \| \overline e^{n+1}\|_\malA \leq {(1+  C \tau)} \left\|e^n\right\|_\malA + \malO \big(\tau^3 + \tau h_x^4+ \tau h_y^4 \big), \quad 0 \leq n \leq N_t-1.
\end{align} 

Finally, the application of discrete Gronwall inequality to \eqref{err:enplus1_end} gives  
\begin{align*}
      \| e^{n+1} \|_\malA \leq C \big(\tau^2+ h_x^4+h_y^4\big),
\end{align*}
if $ \|e^0\|_\malA=\malO \big( h_x^4+ h_y^4 \big)$.
Therefore, Lemma \ref{lem:err:AB_bound} implies the conclusion. \qed
\end{proof}

 \section{Numerical results} \label{sec:test}
In this section, we shall simulate several examples to demonstrate the performance of our proposed schemes, in which the positivity/bound-preserving and mass-conservative properties shall be tested. Meanwhile, the corresponding errors and convergence orders in discrete $L_2$ and $L_{\infty}$-norms shall also be verified. 
In the following tests, we always assume that $N = N_x =N_y$ and thus $h = h_x = h_y$. 
All the numerical experiments are performed using Matlab R2022b on a desktop with the configuration: 12th Gen Intel(R) Core(TM) i5-12400F  @ 2.50GHz and 16 GB RAM. 
	 
\subsection{One-dimensional problems} \label{sec:test:1d}
\begin{example}  (1D viscous Burgers equation)\label{1d:test:burger}
     We consider  the viscous Burgers equation 
   	\begin{align}
     	u_{t}  + f(u)_x= \gamma u_{x x} + S(x,t),  \label{eg:1d:burger}
     \end{align}
    where $ f(u) = {\frac{1}  {2}  u^2}$ and the exact solution is given as one of the following two types:	
\begin{itemize}
	\item Case {1} \cite{xie2008numerical}: 	
    $u(x, t)=\frac{x / t} {1 + \sqrt{t } \exp \left( (x^{2}  -0. 25t) / 4 \gamma t\right) },\quad (x, t) \in [0,3]\times [1,4]$.
     \item Case {2} \cite{acosta2010mollification}: $u(x,t) = 1.5-0.5 \tanh\left( 0.5(x-1.5t)/2\gamma\right), \quad (x, t) \in [-1,3] \times [0,1]$.
 \end{itemize}  
The initial and Dirichlet boundary conditions are taken from the exact solution.
Below, we will illustrate the reliability of our schemes from three aspects: accuracy, bound-preserving property and efficiency. 
\end{example}

To begin with, we verify the accuracy of the HOC-Splitting  scheme  \eqref{1d:scheme:diffu1}--\eqref{1d:scheme:diffu2} and the HOC-Splitting-Substeps scheme \eqref{1d:scheme:diffu1}, \eqref{substep_1}  and \eqref{1d:scheme:diffu2} for Case  {1}. The  discrete  $L_{\infty}$-norm errors and corresponding convergence orders at  $T=4$ for given $\gamma= 5.0\times 10^{-3}$ are presented in Table \ref{Tab:1d:burger:K}. We observe that both schemes without/with substeps bear nearly the same accurate results, and both decrease with fourth-order convergence with respect to (w.r.t.) $h=3/N$.  However,  due to stability limitation, the convergence order for $K=1$ decreases after $N \ge 1024$, but that for the HOC-Splitting-Substeps scheme with $K=10$ remains unchanged. Besides, the HOC-Splitting-Substeps scheme can save about 40\% CPU time compared with the HOC-Splitting scheme without substeps.
\begin{table}[!htbp] 	
	\vskip -0.4cm
	\centering
	\setlength{\abovecaptionskip}{0.1cm} 
	\caption{\small Errors and convergence orders of the HOC-Splitting  scheme   without/with substeps}	\label{Tab:1d:burger:K}
	\begin{tabular}{ccccccc}
		\toprule & \multicolumn{3}{c}{  $K=1$, $\tau$ = 1.0 E-5 } & \multicolumn{3}{c}{  $K=10$, $\tau $ = 1.0 E-4 } \\			\midrule
		$N$    & $L_\infty$-error & Order & CPU time		& $L_\infty$-error & Order  & CPU time\\			\midrule
		180   &  3.6003  E-6  &  ---   & 4.7s	 &  3.6003  E-6  &  ----  & 3.1s   \\    \midrule
		360  &  2.2135  E-7  & 4.02    & 9.3s  &  2.2141  E-7  & 4.02  & 5.1s  \\ 	 \midrule
		540  &   4.3218 E-8  &4.03      & 14.0s  &  4.3323 E-8  & 4.02   & 7.2s  \\      \midrule
		720  &   1.3653 E-8  &4.00      & 17.5s  &  1.3714 E-8  & 4.00   & 9.7s  \\     \bottomrule
	\end{tabular}
\end{table}	
\vskip -0.2cm
	
\begin{figure}[!htbp]\small
	\vspace{-0.4cm}
	\centering  
	\subfigure[ $K = 10, \tau\approx \frac{3h}  { \max |f^{\prime}  |} $ without  limiter ]
	{
		\includegraphics[width=0.45\textwidth]{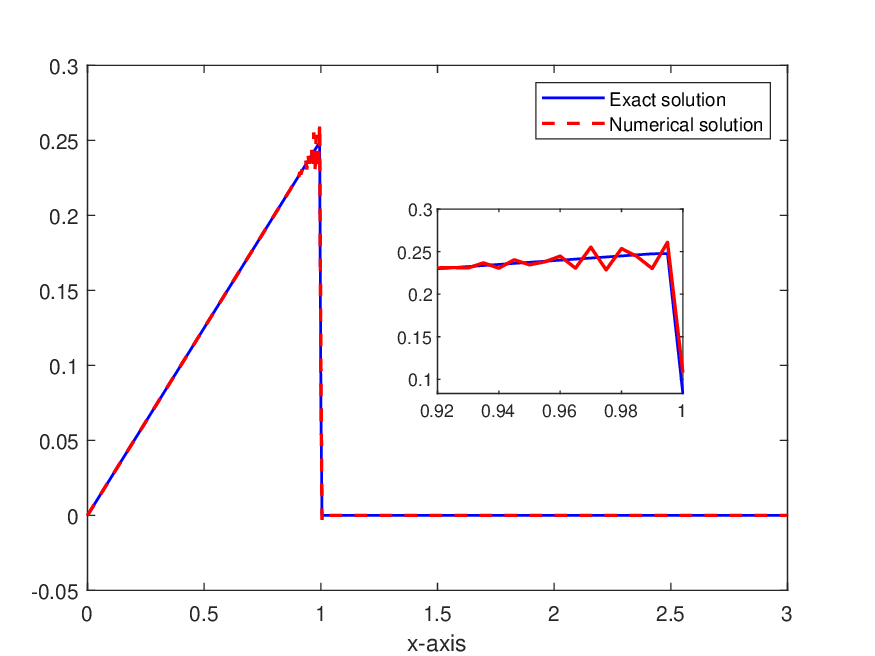}
		\label{Fig:1d:burger:K10non}
	}
	\subfigure[ $K = 1, \tau\approx\frac{h}  { 4 \max |f^{\prime}  |} $ without limiter ]
	{
		\includegraphics[width=0.45\textwidth]{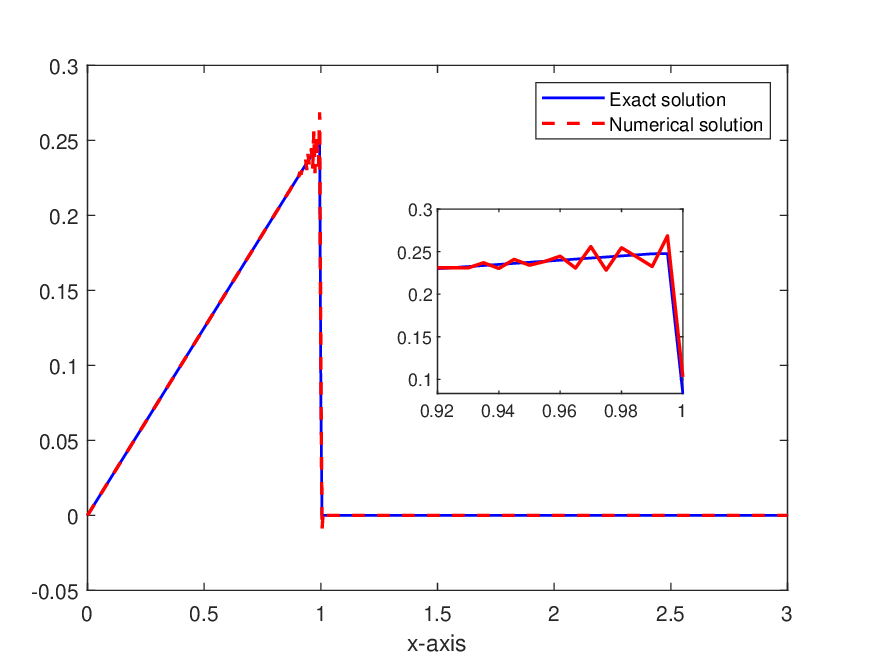}
		\label{Fig:1d:burger:K1non}
	}
	\subfigure[ $K = 10, \tau\approx \frac{3h}  { \max |f^{\prime}  |}   $ with   limiter ]
	{
		\includegraphics[width=0.45\textwidth]{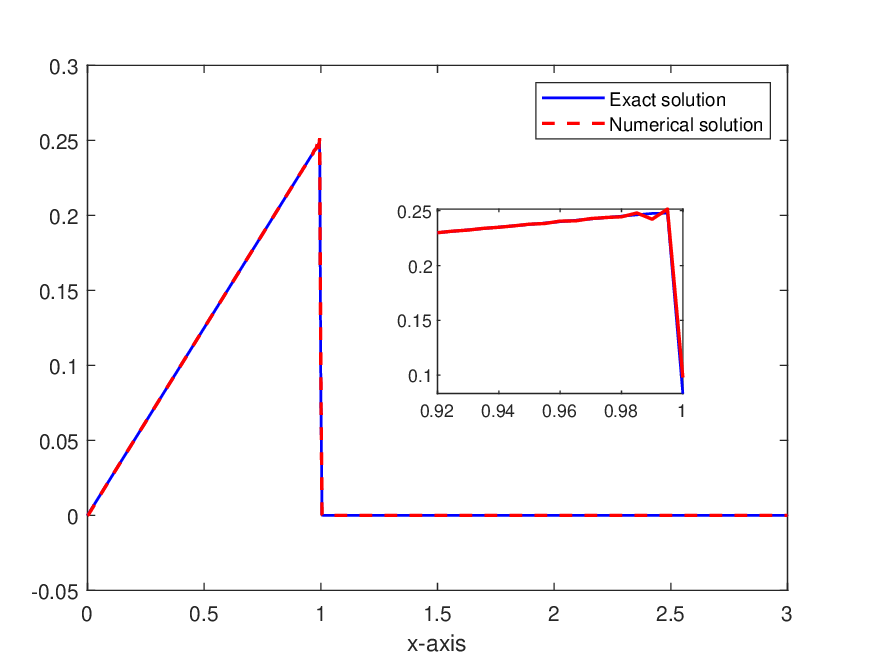}
		\label{Fig:1d:burger:K10lim}
	}
	\subfigure[ $K=1, \tau\approx \frac{h}  {4 \max |f^{\prime}  |}   $ with  limiter ]
	{
		\includegraphics[width=0.45\textwidth]{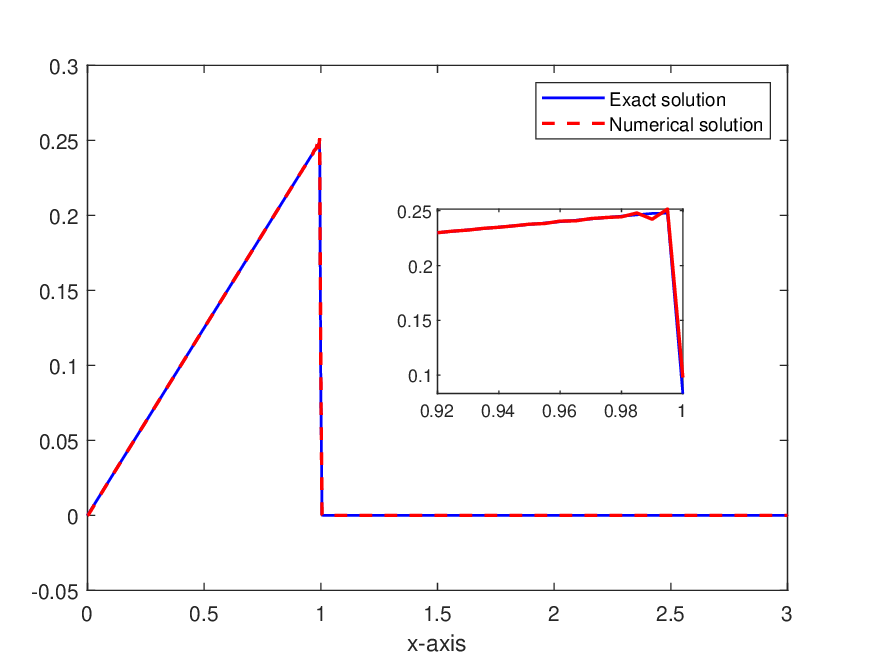}
		\label{Fig:1d:burger:K1lim}
	}
	\setlength{\abovecaptionskip}{0.05cm} 
	\setlength{\belowcaptionskip}{0.0cm}
	\caption{\small  Solution profiles of the HOC-Splitting scheme without/with substeps and TVB limiter}
	\label{Fig:1d:burger:K_lim}
\end{figure}

To further illustrate the effectiveness of the substeps method as well as the TVB limiter, we run the two schemes with small   $\gamma = 1.0\times 10^{-4}$ for  $N=600$, $\tau =6.0\times 10^{-2} \approx \frac{3h}{\max{|f^{\prime}|}} $, where $\max{|f^{\prime}|} = \max{|u|} \approx 0.25$.  
Fig. \ref{Fig:1d:burger:K_lim} indicates that employing the substeps method can further relax the CFL condition, and the TVB limiter can effectively reduce the non-physical oscillations. However, the obtained numerical solution is not bound-preserving, as its minimum value falls below zero. These results are consistent with the conclusions in Remarks \ref{rem:1d:K}--\ref{rem:1d:limter}. 
\begin{figure}  [!htbp]
	\vspace{-0.4cm}
	\centering  
	\subfigure[{HOC-Splitting }  ]
	{ 
		\includegraphics[width=0.45\textwidth]{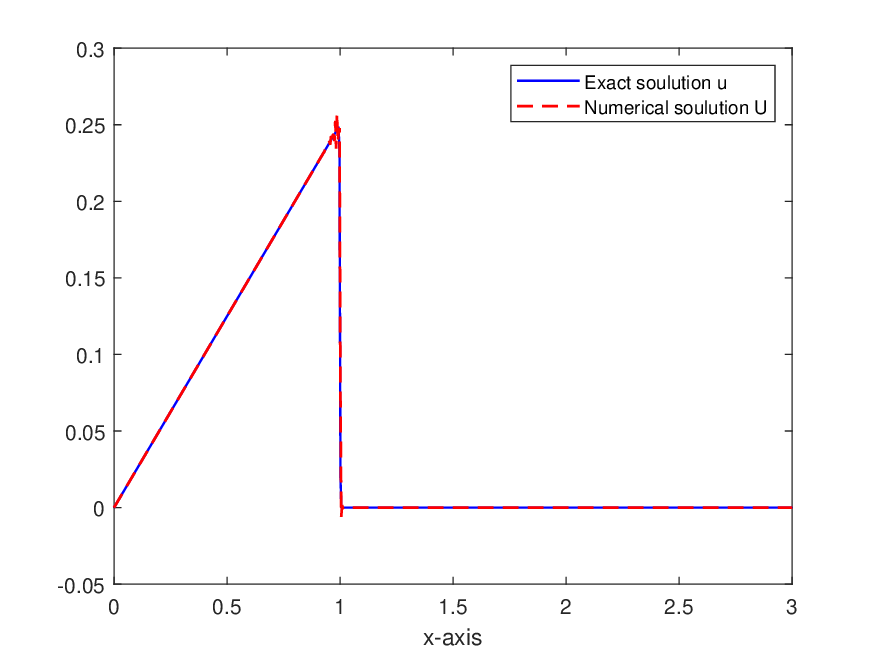} 
          \label{Fig:1d:burger:1e4:non}
		}
	\subfigure[{BP-HOC-Splitting with limiter}  ]
	{
		\includegraphics[width=0.45\textwidth]{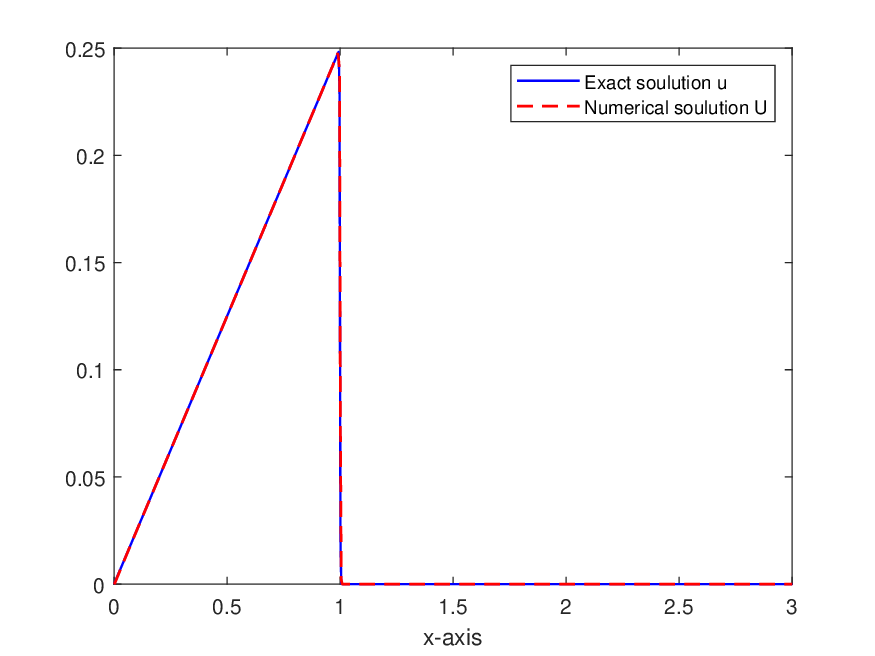}
             \label{Fig:1d:burger:1e4:Lag_lim}
		}
	\setlength{\abovecaptionskip}{0.05cm} 
   \setlength{\belowcaptionskip}{0.0cm}		
	\caption{\small   Comparisons  of  numerical solutions without/with the multiplier $\lambda$ and TVB limiter}	\label{Fig:1d:burger:Lag}
\end{figure}
\begin{table}  [!htbp] 	\small
	\vspace{-0.3cm}
	\centering
	\setlength{\abovecaptionskip}{0.1cm} 
	\caption{\small Errors and convergence orders without/with  Lagrange multiplier $\lambda$}
	\label{Tab:1d:burger:Lag}
	\begin{tabular}{ccccccccc}
		\toprule & \multicolumn{4}{c}{   HOC-Splitting } & \multicolumn{4}{c}{ BP-HOC-Splitting } \\ \midrule 
		$N$  & $L_\infty$-error & Order & $L_2$-error & Order 	& $L_\infty$-error & Order & $L_2$-error & Order\\	\midrule
		90   &  2.8001 E-3  & ---      &  8.3098 E-4  & ---      &  4.2988 E-3  & ---       &  1.0898 E-3  & ---   \\	\midrule
		180  &  2.2571 E-4  & 3.63     &  3.7409 E-5  & 4.47	 &  2.2298 E-4  & 4.26      &  3.7376 E-5  &   4.86\\	\midrule 
		270  &  3.5446 E-5  & 4.56     &  5.8100 E-6  & 4.59	 &  3.5446 E-5  & 4.53      &  5.8100 E-6  & 4.59 \\	\midrule
		360  &  1.0281 E-5  & 4.30     &  1.7181 E-6  & 4.23     & 1.0281 E-5   & 4.30      &  1.7181 E-6  & 4.23\\	\bottomrule
	\end{tabular}
\end{table}	

Next, we test and compare the bound-preserving property of the  HOC-Splitting and BP-HOC-Splitting schemes. Numerical solutions at $T=4$ with $\gamma= 1.0\times 10^{-4}$ and $N = 800, N_t= 600$ are plotted in Fig. \ref{Fig:1d:burger:Lag}.  
We find that the HOC-Splitting scheme appears oscillation and negative solution, see Fig. \ref{Fig:1d:burger:1e4:non}.  As expected, the combination of Lagrange multiplier and TVB limiter can produce nonoscillatory and positivity-preserving numerical solution, as shown in Fig. \ref{Fig:1d:burger:1e4:Lag_lim}. The accuracy of the HOC-Splitting  and BP-HOC-Splitting schemes are also tested, in which the errors and convergence orders  at  $T=4$ are presented in Table \ref{Tab:1d:burger:Lag} for  $\gamma=5\times 10^{-3}$ under $N_t= N^2/3$ (i.e., $\tau = h^2$). As seen, both schemes exhibit second-order temporal accuracy and fourth-order spatial accuracy as stated, indicating that the BP-HOC-Splitting scheme provides a bound-preserving numerical solution without sacrificing accuracy, see also Theorem \ref{th:err:BP-HOC-ADI} for comparison.

\begin{figure}  [!htbp]\small
	\vspace{-0.4cm}
	\centering  
	\subfigure[ $N = 1000, \tau\approx \frac{h}  {3 \max |f^{\prime}  |}   $ without multiplier]
	{
		\includegraphics[width=0.45\textwidth]{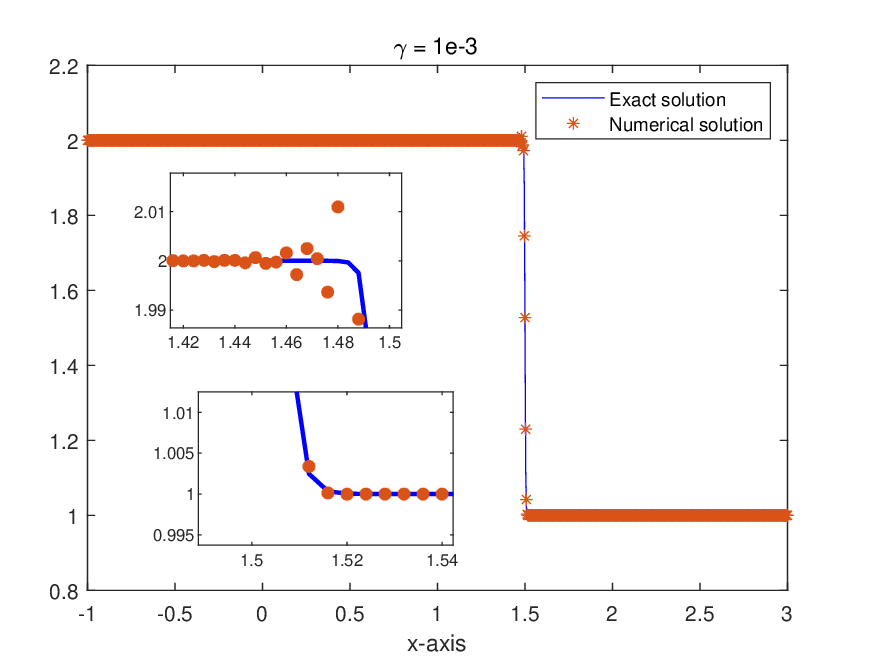}
	}
	\subfigure[ $N = 1000, \tau\approx \frac{h}  {3 \max |f^{\prime}  |}   $ with   multiplier]
	{
		\includegraphics[width=0.45\textwidth]{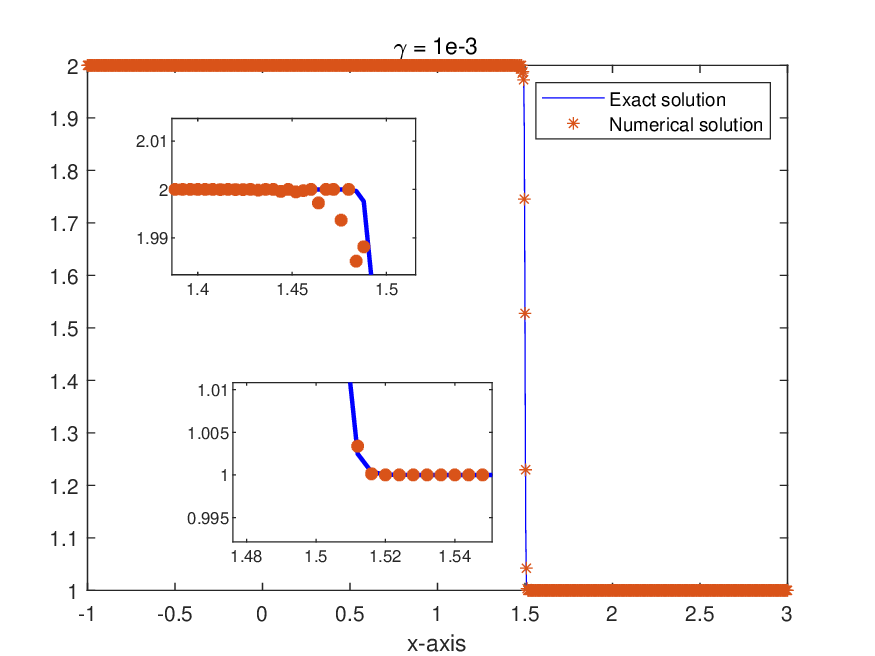}
	}
	\setlength{\abovecaptionskip}{0.05cm} 
	\setlength{\belowcaptionskip}{0.0cm}
	\caption{\small Comparison of solutions obtained  without/with  the multiplier $\lambda$}
	\label{fig:1d:bur:Diriclet}
\end{figure}

In addition, we also test the bound-preserving property of the BP-HOC-Splitting scheme for Case  {2},  in which non-homogeneous boundary condition is considered. Fig. \ref{fig:1d:bur:Diriclet}  shows the exact and approximate solutions at $t =1$  without/with the multiplier $\lambda$  for $\gamma=1.0\times 10^{-3} $.  Note that the exact solution should always fall within the range $[1, 2]$, and meanwhile, a large gradient appears in the solution and is moving w.r.t. time, which requires our scheme to be bound-preserving and can better characterize the large-gradient solution. It is observed that the numerical solution computed by the BP-HOC-Splitting scheme coincides with the exact solution and maintains the bound-preserving property. In contrast, the numerical solution without using the multiplier exceeds the upper and/or lower bounds in the neighborhood of the large gradient.

\begin{example}\label{1d:test:FP} (1D Fokker-Planck equation) In this example, we consider the simulation of Fokker-Planck equation \cite{cheng2022new2}:
	\begin{align}
		\partial_t u=\partial_x \left( x u\left( 1-u\right)+\partial_x u\right),\quad x\in [-2 \pi, 2 \pi], \label{fp1}
	\end{align}
under periodic boundary condition and initial condition $u(x, 0)= \exp ({-x^2/0.4})$. 
In real applications, the equation is usually used to model the relaxation of fermion and boson gases taking on the form described in \cite{carrillo20081d}, and it can also be interpreted as a gradient flow:
	\begin{align*}
		\partial_t u=\partial_x \Big( u\ ( 1-u)\partial_x \frac{\delta \malE}{\delta u}\Big) 
	\end{align*}
	with $\malE=\malE(u)$ being the entropy functional
	\begin{align}
		\malE(u)=\int_{\Omega} \Big( \frac{x^2}{2} u+u \log  u+(1-u) \log ( 1-u) \Big)  d x. \label{fp_E}
	\end{align}
 It is also important to note that the solution of \eqref{fp1} is expected to always take values in $[0,1]$, since negative values would render the energy functional \eqref{fp_E} meaningless due to the terms $\log u$ and $\log \left( 1-u\right)$ \cite{cheng2022new2}. Moreover, this equation with periodic boundary condition always conserves the mass, i.e., $\frac{d}{dt}\int_{\Omega} u dx=0$. Therefore, hopefully it requires the numerical solution to be bound-preserving and mass-conservative. 
 \end{example}
 
\begin{figure}  [!htbp]\small
		\vspace{-0.5cm}
	\centering  
	\subfigure 
	{
		\includegraphics[width=0.45\textwidth]{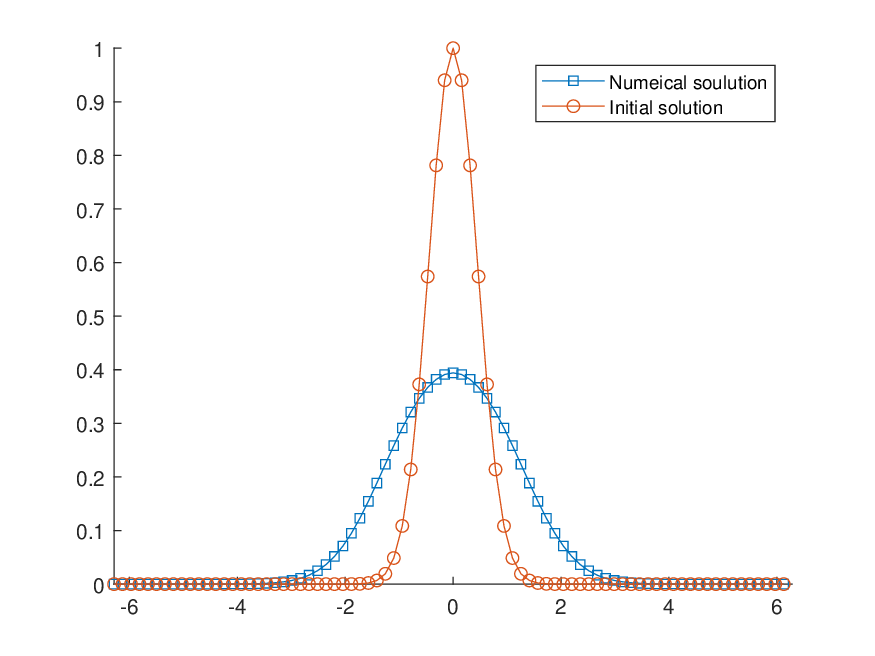}
	}
	\subfigure 
	{
		\includegraphics[width=0.45\textwidth]{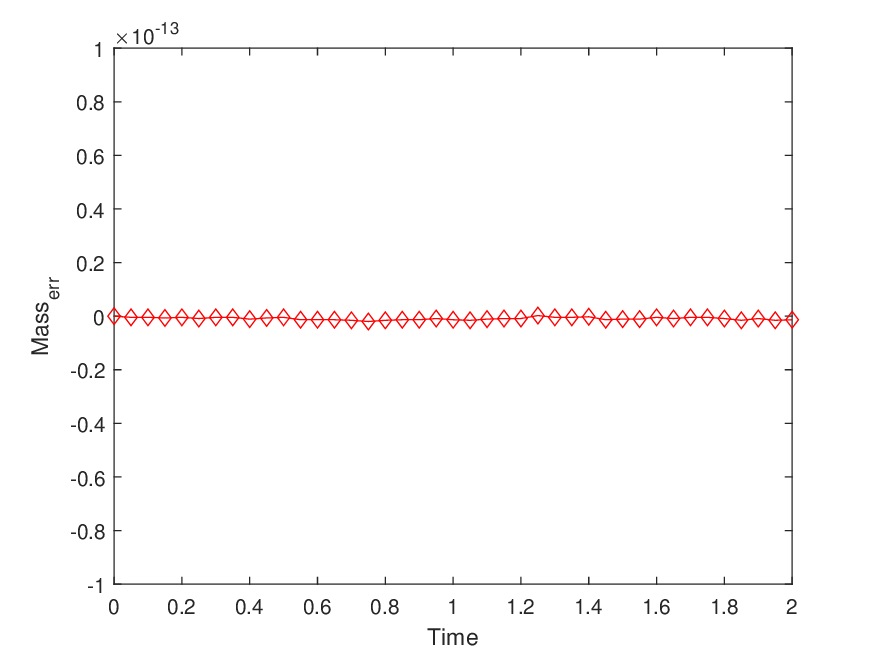}
	}
	\setlength{\abovecaptionskip}{0.05cm} 
	\setlength{\belowcaptionskip}{0.0cm}
	\caption{\small   Numerical solution (left) and  evolution of mass errors w.r.t. time (right)}
	\label{1d:FP:compare}
\end{figure}

\begin{figure}  [!htbp]\small
	\vspace{-1.0cm}
	\centering
	\includegraphics[width=0.45\textwidth]{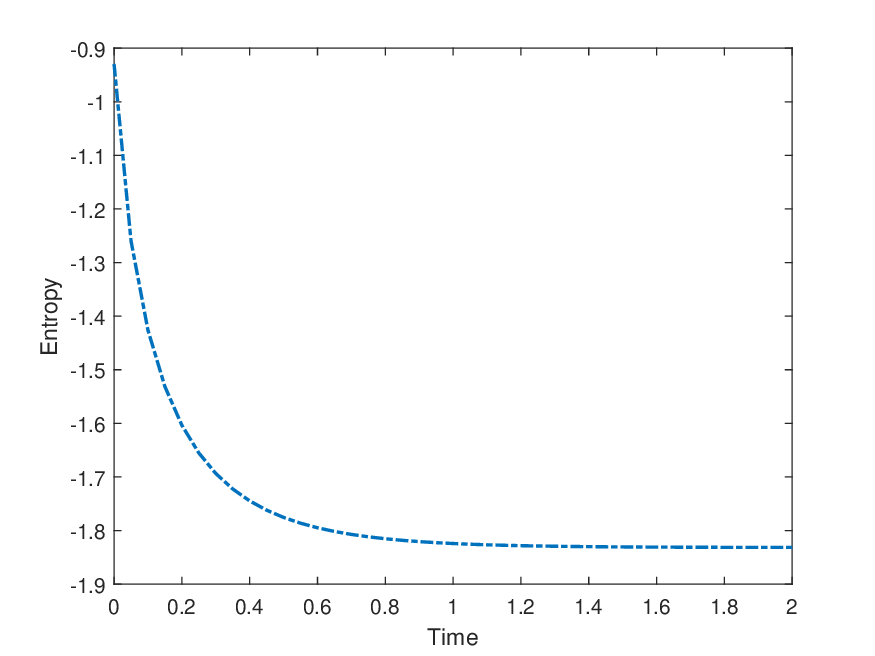}
		\setlength{\abovecaptionskip}{0.05cm} 
	\setlength{\belowcaptionskip}{0.0cm}
	\caption{\small   Evolution of the discrete entropy}
	\label{Fig:1d:FP:Energy}
\end{figure}

We now run the BP-MC-HOC-Splitting scheme to check these properties. Numerical solution profile at time $T=2$ and  the evolution of  mass errors w.r.t. time for $N=80, N_t = 40$ are plotted in Fig. \ref{1d:FP:compare}. It can be seen that  the scheme indeed generates numerical solutions remain within the interval $[0,1]$ and satisfies the mass-conservative relation. 
In addition, it is well-known that the Fokker-Planck  equation \eqref{fp1}  satisfies the entropy dissipation law, i.e., $\frac{d }{dt} \malE(u)\leq 0$.  We depict the discrete version entropy for the BP-MC-HOC-Splitting scheme in Fig. \ref{Fig:1d:FP:Energy}.  It is observed that the discrete entropy is dissipative, which strongly supports the effectiveness of the scheme.

In the simulation, we also examine the accuracy of the BP-MC-HOC-Splitting scheme. Since there is no available analytical solution, we compute the Cauchy error \cite{Cauchyerror} via the numerical solution at $T=0.5$, which is displayed in Tables \ref{Tab:1d:FP:Space}--\ref{Tab:1d:FP:Time}. 
 It clearly  shows that the proposed scheme is fourth-order accurate in space and second-order accurate in time under the discrete $L_\infty$-norm, while the discrete $L_2$-norm temporal convergence orders can reach 2.5, indicating a possible superconvergent result.
\begin{table}  [!htbp] 	\small
		\vspace{-0.5cm}
	\centering
	\setlength{\abovecaptionskip}{0.1cm} 
	\caption{\small  Spatial errors of the BP-MC-HOC-Splitting scheme with $\tau=10^{-4}$}
	\label{Tab:1d:FP:Space}
	\begin{tabular}{lllll }
		\toprule $ N$             & 50            & 100           & 200          & 400 \\
		\midrule $L_\infty$-error & 2.6451  E-5   &  1.5731  E-6  & 9.7181  E-8  &  6.0576  E-9 \\
		Order                     & ---           &  4.07         &  4.01        & 4.00  \\
		\midrule $L_2$-error      & 2.6959  E-5   &  1.6264  E-6  & 1.0091  E-7  &  6.2964  E-9  \\
		Order                     & ---           & 4.05          &  4.01        & 4.00 \\		
		\bottomrule
	\end{tabular}	
\end{table}	
\begin{table}  [!htbp] 	\small
	\centering
	\setlength{\abovecaptionskip}{0.1cm} 
	\caption{\small  Temporal errors of the BP-MC-HOC-Splitting scheme with $h=0.025$}
	\label{Tab:1d:FP:Time}
	\begin{tabular}{lllll}
		\toprule $N_t$            & 100          & 200         & 400          & 800 \\
		\midrule $L_\infty$-error &  1.0858 E-6  & 2.7156 E-7  &  6.7906 E-8  & 1.6978 E-8 \\
		Order           & ---          &  2.00       & 2.00         & 2.00  \\
		\midrule $L_2$-error      &  6.6077 E-7  & 1.1679 E-7  &  2.0645 E-8  & 3.6495 E-9  \\
		Order           & ---          & 2.50        & 2.50         & 2.50  \\
		\bottomrule
	\end{tabular}	
\end{table}	

\subsection{  Two-dimensional problems} \label{sec:test:2d}
\begin{example}\label{2d:test:burgers}  (2D viscous  Burgers equation)	Consider the  2D Burgers-type equation with periodic boundary condition \cite{elton1996comparisons}
			\begin{equation}
				u_{t}  +f(u)_{x}  + g(u)_{y}   = \gamma (u_{x x}  +u_{y y}  ) + S(x,y,t), \quad
				(x, y, t) \in[-1, 2]^2 \times[0, 0. 6], \notag
			\end{equation}
			where  $f(u)= g(u) = \frac{1}  {2}   u^2$.  Given the source term such that  the  exact solution  is
			\begin{align}   \label{2d:eg:burger:sol}
				u(x, y, t)=\frac{\sigma^{2}  }  { \sigma^{2}  +2\gamma t}   \exp\left( -\frac{({x}  -t - x_0)^{2}  +({y} -t - y_0)^{2}  }  {2 \sigma^{2}  +4 \gamma t}  \right),
			\end{align}
with $\sigma= 0.07$ and $x_0=y_0=0.5$.	
\end{example}

Note that the equation is  mass conservative in the sense of \eqref{model:mass}, and the solution to this problem is positive. Therefore,  the positivity and mass conservation of the numerical solution will  be tested in this example.	 Furthermore, in order to better illustrate the efficiency of the proposed schemes, we also implement the following  second-order implicit-explicit HOC scheme (BDF2-IMEX-HOC):
       \begin{equation}\label{2d:BDF2-HOC:n}
    	  \begin{aligned}
	    	& \frac{ 3 u_{i,j}^{n+1}-  4u_{i,j}^n  + u_{i,j}^{n-1}  }{2\tau } - \gamma ( \malA_{x}^{-1} \delta_ x^2 + \malA_{y}^{-1} \delta_ y^2)  u_{i,j}^{n+1} \\
	    	& \qquad + \mathcal{B}_x^{-1} D_{\hat{x}} f(2u_{i,j}^{n}-u_{i,j}^{n-1} ) +   \mathcal{B}_y^{-1} D_{\hat{y}} g(2u_{i,j}^{n}-u_{i,j}^{n-1} ) =S_{i,j}^{n+1}, ~~ n \ge 1,
	     \end{aligned}
      \end{equation}
and $u^1$  is computed via the first-order implicit-explicit HOC scheme  
\begin{align}\label{2d:BDF2-HOC:1}
	\begin{aligned}
		& \frac{  {u}_{i,j}^{1}-  u_{i,j}^0  }{\tau } - \gamma ( \malA_{x}^{-1} \delta_ x^2 + \malA_{y}^{-1} \delta_ y^2)  {u}_{i,j}^{1} + \mathcal{B}_x^{-1} D_{\hat{x}} f(u_{i,j}^{0} ) +   \mathcal{B}_y^{-1} D_{\hat{y}} g(u_{i,j}^{0} )  =S_{i,j}^{1}.
	\end{aligned}
\end{align}

\begin{table}  [!htbp] \small
		\vspace{-0.6cm}
	\centering
	\setlength{\abovecaptionskip}{0.1cm} 
	\caption{\small Comparisons of  BDF2-IMEX-HOC and HOC-ADI-Splitting schemes}
	\label{tab:2D:compare_order}
	\begin{tabular}{lcccccc}
		\toprule 
		& \multicolumn{3}{c}{ BDF2-IMEX-HOC  } & \multicolumn{3}{c}{ HOC-ADI-Splitting } \\		\midrule 
		$N$  & $L_\infty$-error & Order  & CPU time 	& $L_\infty$-error & Order  & CPU time \\		\midrule
		90  & 3.8444 E-4  & ---   &  12.9s  	& 4.1895 E-4  & ---   &  0.9s  \\
		180 & 1.9813 E-5  & 4.27  &  4m 6s		& 2.2033 E-5  & 4.24  &  15s \\
		270 & 3.7994 E-6  & 4.07  &  25m 13s  	& 4.2384 E-6  & 4.06  &  1m 7s \\			
		360 & 1.1902 E-6  & 4.03  &  86m 20s	 & 1.3291 E-6  & 4.03  &  6m 18s\\
		\bottomrule
	\end{tabular}
\end{table}	
In the following, we first test the accuracy of  the  BDF2-IMEX-HOC scheme  \eqref{2d:BDF2-HOC:n}--\eqref{2d:BDF2-HOC:1} and     HOC-ADI-Splitting scheme \eqref{HOC-ADI:e1a}--\eqref{HOC-ADI:e4b}  with $\gamma = 5.0\times 10^{-3}$.   The corresponding numerical results are listed in Table \ref{tab:2D:compare_order} for   $N^2=15 N_t $  (i.e., $\tau = h^2$), which actually indicates that  both methods generate almost the same accurate numerical errors and have fourth-order accuracy in space and second-order accuracy in time. However, we observe that the proposed splitting method has significantly improved computational efficiency. For example, the implementation of the BDF2-IMEX-HOC scheme takes about 14 times  CPU time of the proposed HOC-ADI-Splitting scheme  when $ N =  360$.

\begin{figure} [!htbp] \small
	\vspace{-0.5cm}
	\centering
	\subfigure[ {   Exact  solu.} ]
	{
		\includegraphics[width=0.45\textwidth]{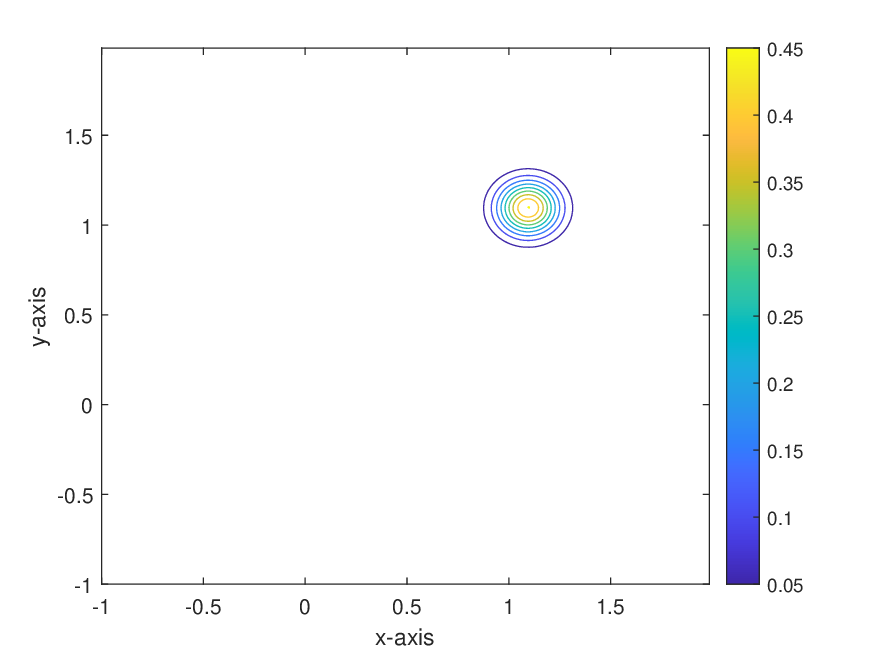}
	}
	\subfigure[ {  HOC-ADI-Splitting solu.}  ]
	{
		\includegraphics[width=0.45\textwidth]{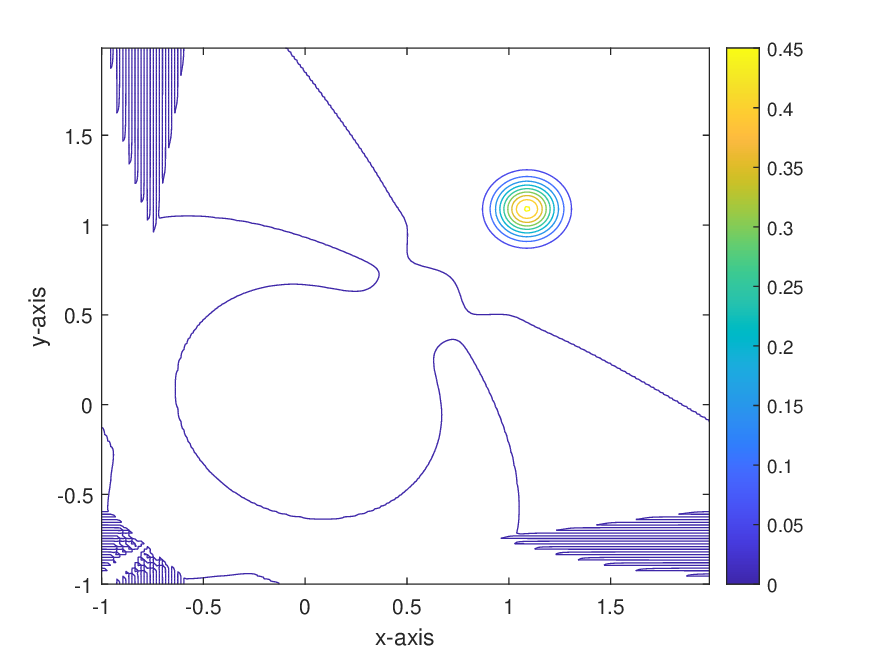}
	}
	\subfigure[ { BP-HOC-ADI-Splitting solu.}  ]
	{
		\label{Fig. sub. 1}
		\includegraphics[width=0.45\textwidth]{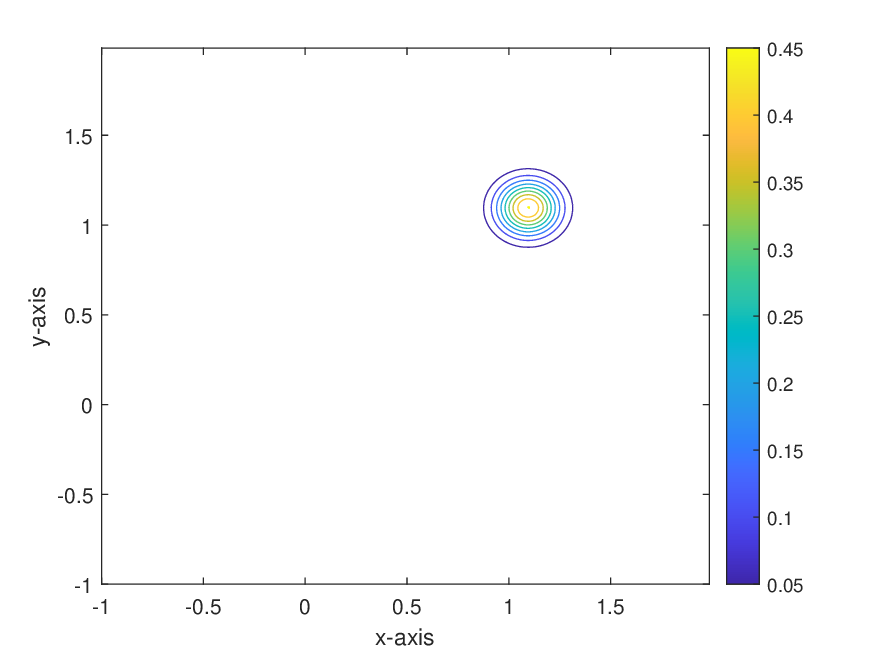}
	}
	\subfigure[ { BP-MC-HOC-ADI-Splitting solu.}  ]
	{
		\label{Fig. sub. 2}
		\includegraphics[width=0.45\textwidth]{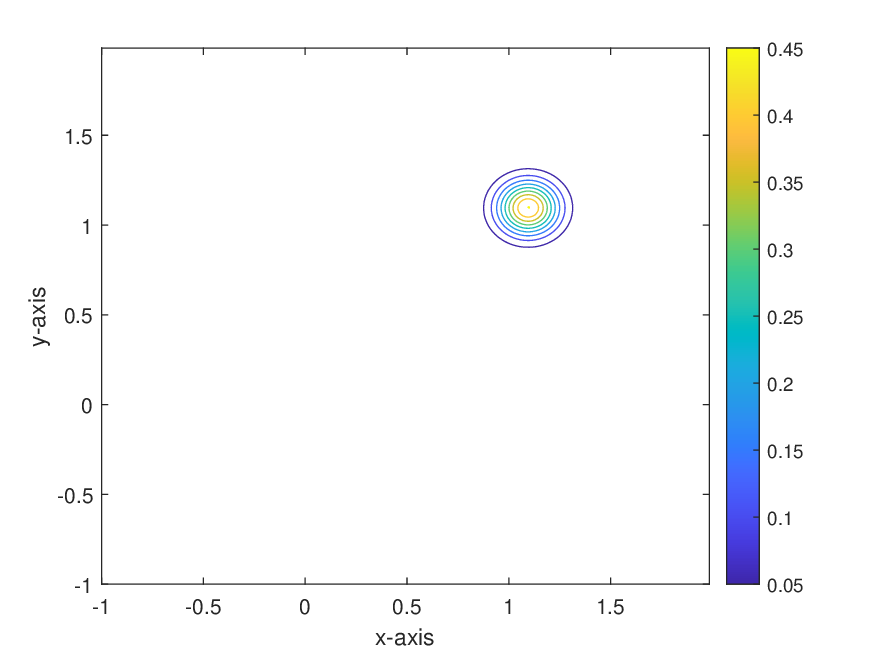}
	}	
	\setlength{\abovecaptionskip}{0.05cm} 
	\setlength{\belowcaptionskip}{0.0cm}
	\caption{\small   Contour plots of computed solutions  without/with Lagrange multipliers}
	\label{fig:2d:burger:u:contour}
\end{figure}

\begin{figure}[!hptb] \small
 	\vspace{-0.4cm}
		\centering
		\subfigure[ { HOC-ADI-Splitting}   ]
	{
		\includegraphics[width=0.31\textwidth]{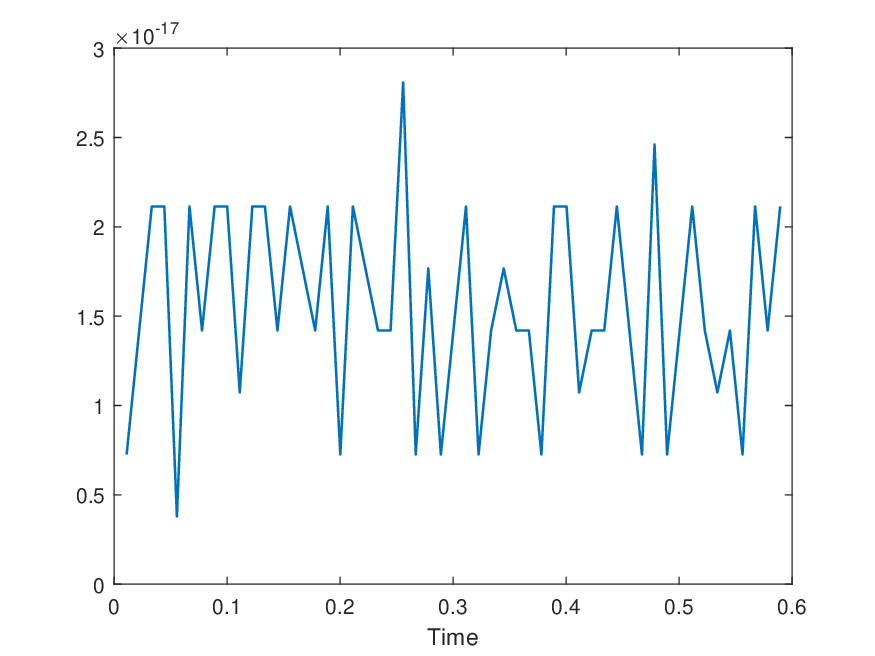}
	}
	    \subfigure[ { BP-HOC-ADI-Splitting}  ]
	{
     	\includegraphics[width=0.31\textwidth]{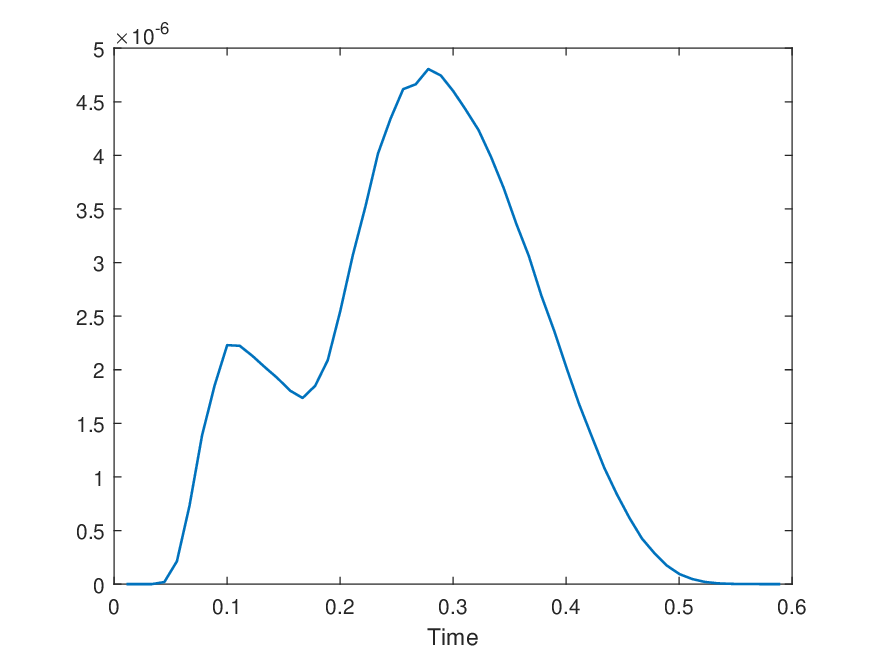}
	}
		\subfigure[ { BP-MC-HOC-ADI-Splitting}  ]
	{
		\includegraphics[width=0.31\textwidth]{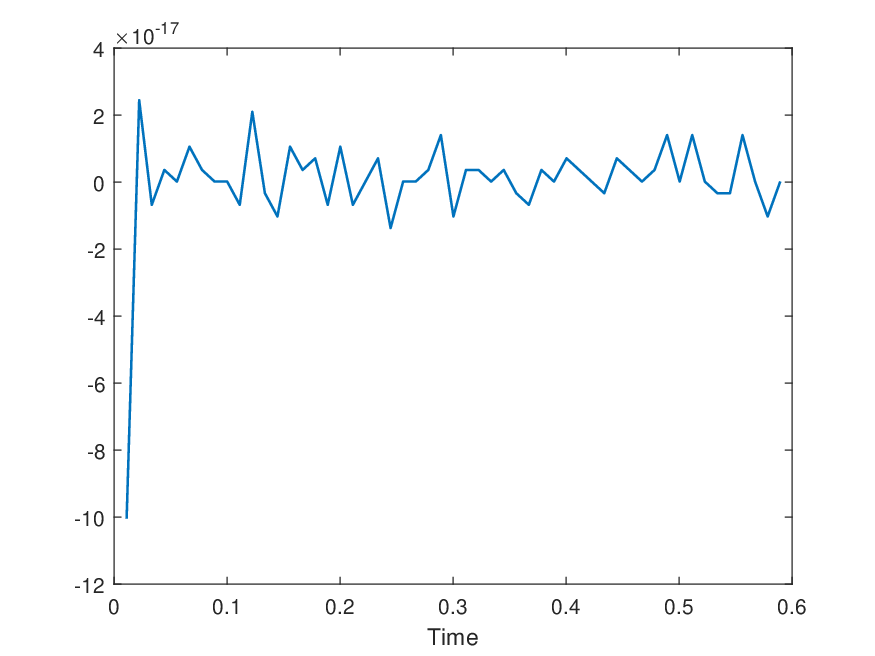}
	}
		\setlength{\abovecaptionskip}{0.05cm} 
	\setlength{\belowcaptionskip}{0.0cm}
	\caption{\small   The numerical mass errors without/with Lagrange multipliers}
	\label{fig:2d:burger:u:mass}
	\end{figure}

Next, we verify the bound-preserving and mass-conservation properties of the  HOC-ADI-Splitting schemes  without/with Lagrange multipliers. Fig. \ref{fig:2d:burger:u:contour} shows the comparisons of numerical solutions obtained by the HOC-ADI-Splitting  scheme \eqref{HOC-ADI:e1a}--\eqref{HOC-ADI:e4b}, the BP-HOC-ADI-Splitting scheme \eqref{HOC-ADI:BP:e1a}--\eqref{HOC-ADI:BP:e2} and the BP-MC-HOC-ADI-Splitting scheme  \eqref{HOC-ADI:Mass:e1}--\eqref{HOC-ADI:Mass:correct}, where $\gamma = 5.0\times 10^{-3}$, $N= 200$ and $\tau = 5.6  \times 10^{-3} \approx\frac{h}{3\max \{|f^{\prime}|,|g^{\prime}|\}}$. 
It shows that polluted negative values appear if the bound-preserving technique is not applied, while the BP-MC-HOC-ADI-Splitting scheme can remedy the negative values in a conservative way. Moreover, though the BP-HOC-ADI-Splitting scheme can preserve positivity, it is not mass-conservative due to the cut-off used in \eqref{HOC-ADI:BP:e2}, see Fig. \ref{fig:2d:burger:u:mass}.

Finally, we also test the convergence rates of the BP-HOC-ADI-Splitting scheme \eqref{HOC-ADI:BP:e1a}--\eqref{HOC-ADI:BP:e2} and BP-MC-HOC-ADI-Splitting scheme \eqref{HOC-ADI:Mass:e1}--\eqref{HOC-ADI:Mass:correct}. The errors in the discrete $L_{\infty}$ and $L_2$ norms at $T=0.6$  with $\tau=h^2$, together with the corresponding convergence orders, are presented in Table \ref{2d:burger:triangle}. It can be found that both of the schemes hold fourth order in space and second order in time, which agrees with the theoretical results in Sect. \ref{sec:ErrEst}. It further shows that both the  BP-HOC-ADI-Splitting scheme  and  BP-MC-HOC-ADI-Splitting scheme  can achieve the corresponding structure-preserving properties without destroying the accuracy.  
\begin{table}  [!htbp] 	 \small
	\centering
	\setlength{\abovecaptionskip}{0.1cm} 
	\caption{\small Errors of the BP-HOC-ADI-Splitting (I) and BP-MC-HOC-ADI-Splitting (II) schemes}	\label{2d:burger:triangle}
	\begin{tabular}{l|lllll}
		\toprule   &  $N$          
    &  60          &  120        & 240         & 360  \\
	\midrule 	   
 I
    &   $L_\infty$-error 
    &  2.6914 E-3  & 1.2300 E-4  & 7.0701 E-6  & 1.3776 E-6 \\
		    &  Order                 
    & ---          &  4.45       & 4.12        & 4.03  \\
		  &    $L_2$-error      
    &  3.8533 E-4  & 1.9302 E-5  & 1.1231 E-6  &  2.1913 E-7  \\
		    &  Order               
    & ---          & 4.32        & 4.10         & 4.03  \\
		\midrule   
  II
  &    $L_\infty$-error 
    &  2.5986 E-3  & 1.1965 E-4  & 6.9191 E-6  & 1.3521 E-6 \\
		      &  Order             
    & ---          &  4.44       &  4.11         & 4.03  \\
		   &    $L_2$-error      
    &  4.0751 E-4  & 1.9827 E-5  &   1.1412 E-6  &  2.2165 E-7  \\
		      &  Order               
    & ---         & 4.38         & 4.12         & 4.04  \\
		\bottomrule
	\end{tabular}	
\end{table}

  \begin{figure}  [!ht] \small
 	\vspace{-0.3cm}
 	\centering 
 	\subfigure[{ HOC-ADI-Splitting}  ]
 	{
 		\includegraphics[width=0.45\textwidth]{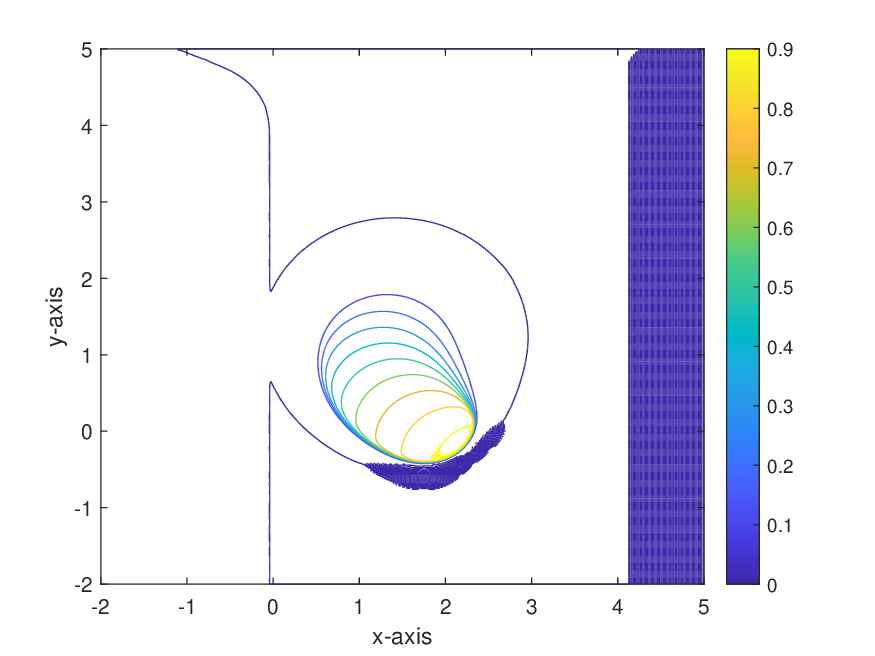}
 	}
 	\subfigure[{  BP-MC-HOC-ADI-Splitting }  ]
 	{	
 		\includegraphics[width=0.45\textwidth]{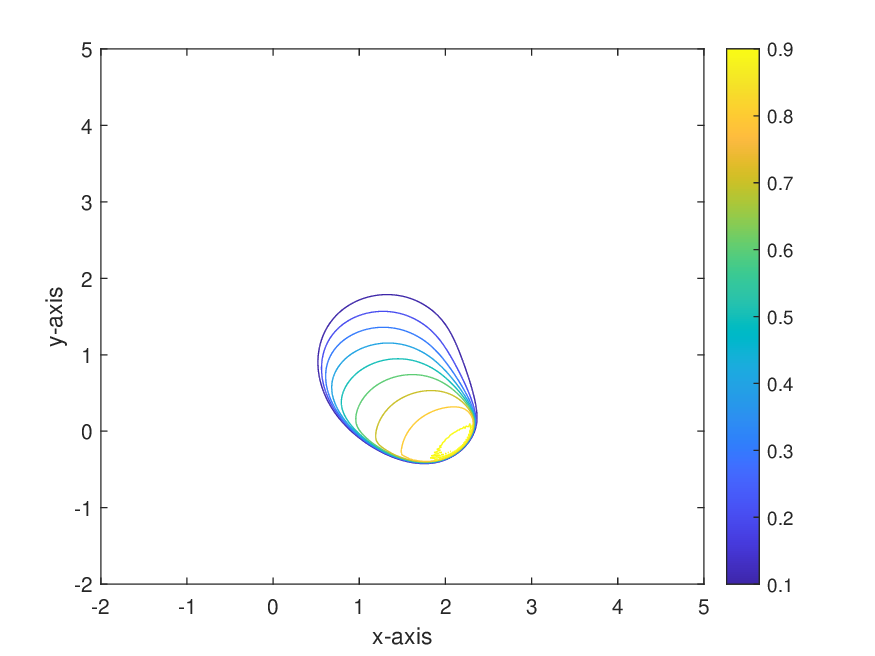}
 	}     
 	\subfigure[{ HOC-ADI-Splitting  }  ]
 	{
 		\includegraphics[width=0.45\textwidth]{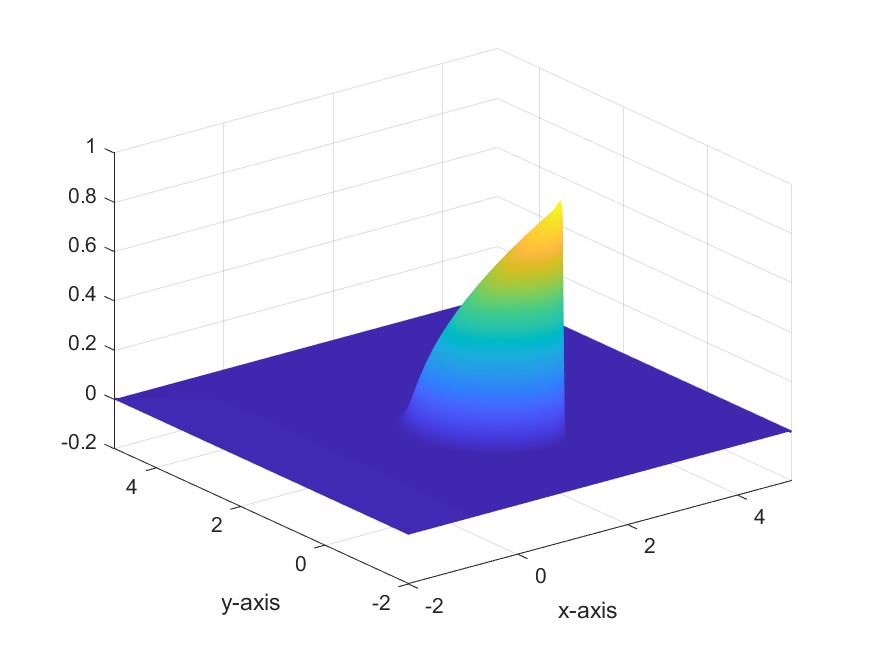}
 	}
 	\subfigure[{ BP-MC-HOC-ADI-Splitting}  ]
 	{
 		\includegraphics[width=0.45\textwidth]{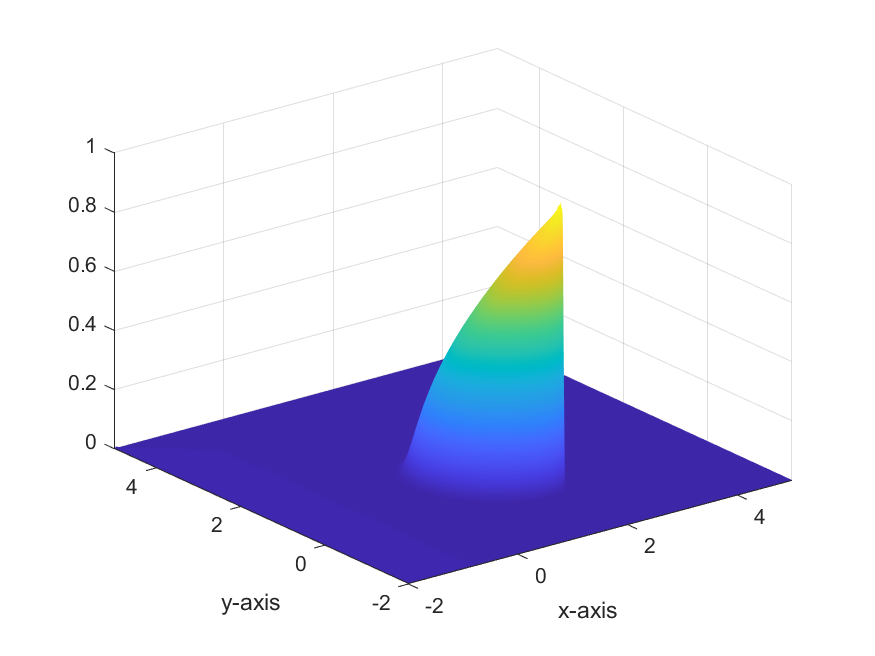}
 	}
 	\setlength{\abovecaptionskip}{0.05cm} 
 	\setlength{\belowcaptionskip}{0.0cm}
 	\caption{\small   Contour plots and three-dimensional views of numerical solutions without/with the multipliers}
 	\label{fig2d_BL_ref}
 \end{figure}
\begin{example}\label{2d:test:BL}   (Buckley-leverett type equation)  
 Consider the 2D nonlinear convection diffusion equation \cite{karlsen1997operator}
  \begin{equation*}
	u_{t}  +(u+(u-0. 25)^{3}  )_{x}  -(u+u^{2}  )_{y}  =\gamma (u_{x x}  +u_{y y}  ), \\
		\quad
	(x, y, t) \in[-2, 5]^2 \times[0, 1],   
 \end{equation*}
with initial data given by
\begin{align*}
 u(x,y,0) = \begin{cases}  1, &   (x-0. 25)^{2}  +(y-2. 25)^{2}  <0. 5, \\ 0, & \text { otherwise.  }  \end{cases}
\end{align*}
The boundary values are set to zero, which is consistent with the initial data. 
\end{example}
 
Although the initial solution is discontinuous, the solution itself is bound-preserving and $\min(u,0)=0$ \cite{karlsen1997operator}. We compute the numerical solution using the HOC-ADI-Splitting  scheme \eqref{HOC-ADI:e1a}--\eqref{HOC-ADI:e4b} and the BP-HOC-ADI-Splitting scheme \eqref{HOC-ADI:BP:e1a}--\eqref{HOC-ADI:BP:e2}, respectively. To improve the computational efficiency, we have adopted the substeps method \eqref{substep_1} with $K=10$ in all simulations. Comparisons of the contour plots and three-dimensional views of numerical solutions obtained by the HOC-ADI-Splitting  scheme and the BP-HOC-ADI-Splitting scheme at $T=1$ are plotted in Fig. \ref{fig2d_BL_ref} for $N = N_t = 500$ and $\gamma = 5.0\times 10 ^{-3}$. It is evident that the proposed scheme, which incorporates the Lagrange multiplier $\lambda$, can effectively preserve the bounds. The results obtained with the multiplier are consistent with the conclusions reported in \cite{karlsen1997operator}, thus validating the effectiveness of our numerical scheme.	

\begin{figure}[!ht]
	\vspace{-0.6cm}
	\centering
	\includegraphics[width=0.45\textwidth]{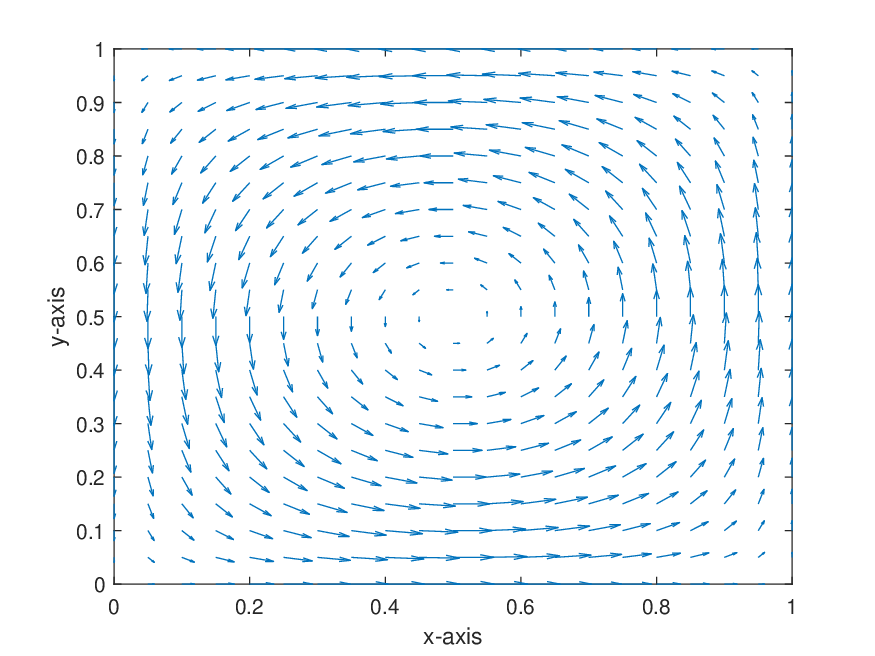}
	\includegraphics[width=0.45\textwidth]{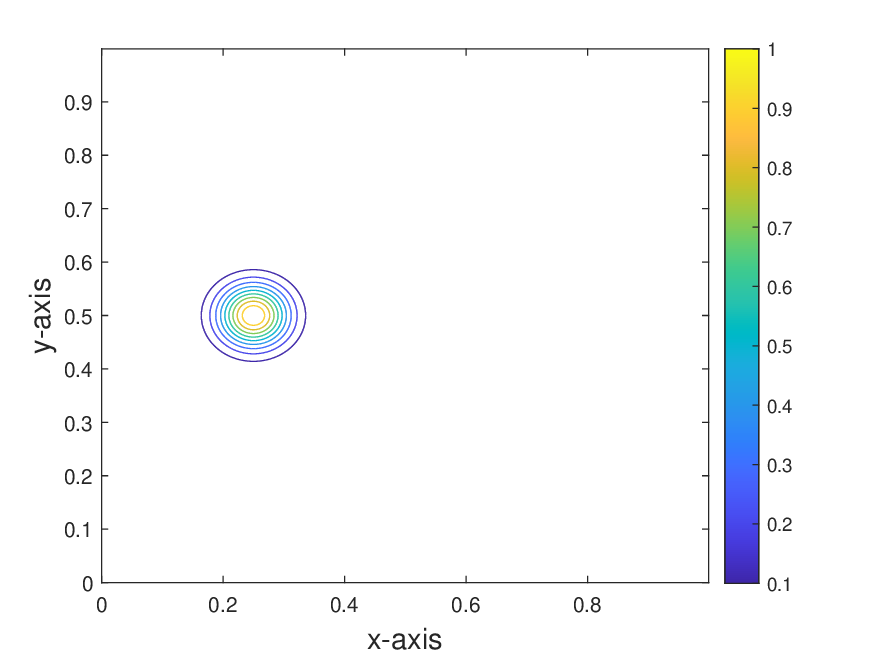}
	\setlength{\abovecaptionskip}{0.05cm} 
	\setlength{\belowcaptionskip}{0.0cm}
	\caption{\small  Contour plots of  Vortex velocity field (left) and  initial Gaussian hump (right)}
	\label{fig:roate:initial}
\end{figure}

\begin{example}\label{2d:test:Gaussian}   (Gaussian hump in vortex shear) 
Consider the transport of a Gaussian hump governed by the linear 2D variable-coefficient convection diffusion problem \cite{fu2017time}
\begin{equation*}
     u_t   +  (c_1 u)_ x  + (c_2 u)_y  = \gamma  (u_{xx}+  u_{yy}),  \quad (x, y, t) \in[0, 1]^2 \times[0, 2],
\end{equation*}
with a divergence-free velocity field $(c_1,c_2)  = (\psi _y, -\psi _x)$ such that $\psi (x, y) = \frac{1}{\sqrt{2}}\sin(\pi x)\sin(\pi y)$, and the initial distribution of the Gaussian hump is given as follows:
\begin{equation*}
u^o(x, y)= \exp \Big(-\frac{( {x} -x_0)^2+( {y} -y_0)^2}{2 \sigma^2 }\Big),
\end{equation*}
where the center of the initial Gaussian hump is specified as $(x_0, y_0) = (0.25, 0.5)$,  with $\sigma^2=1.6 \times 10^{-3}$ as the standard deviation and $\gamma = 5 \times 10^{-5}$.  The contour plots of the velocity field and the initial Gaussian hump are illustrated in Fig. \ref{fig:roate:initial}. 
\end{example}
\begin{figure}  [!htbp]
	\vskip -0.4cm
	\centering   
	\subfigure[$t=0.5$]
	{
		\includegraphics[width=0.3\textwidth]{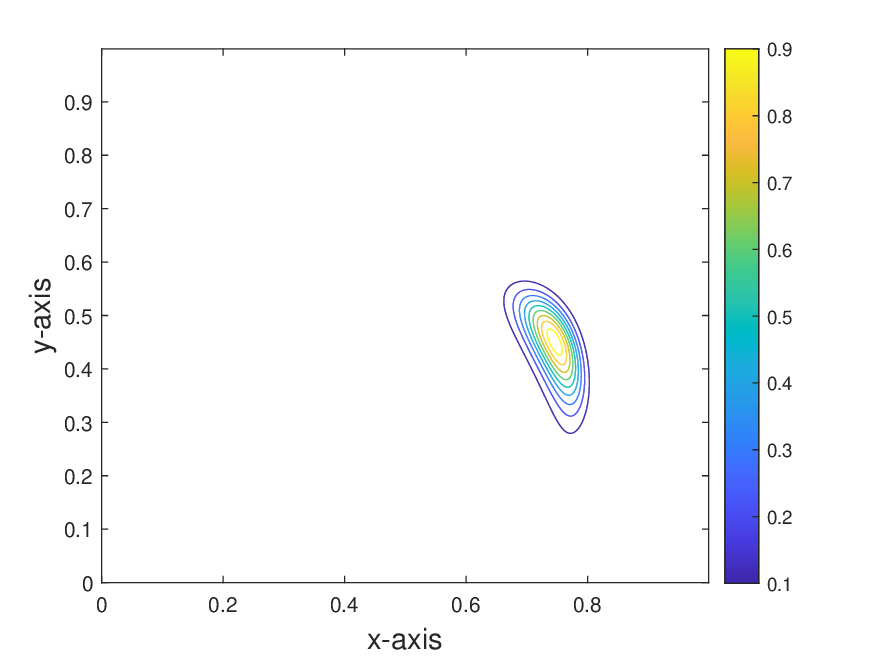}
	}
	\subfigure[$t=1.5$]
	{
		\includegraphics[width=0.3\textwidth]{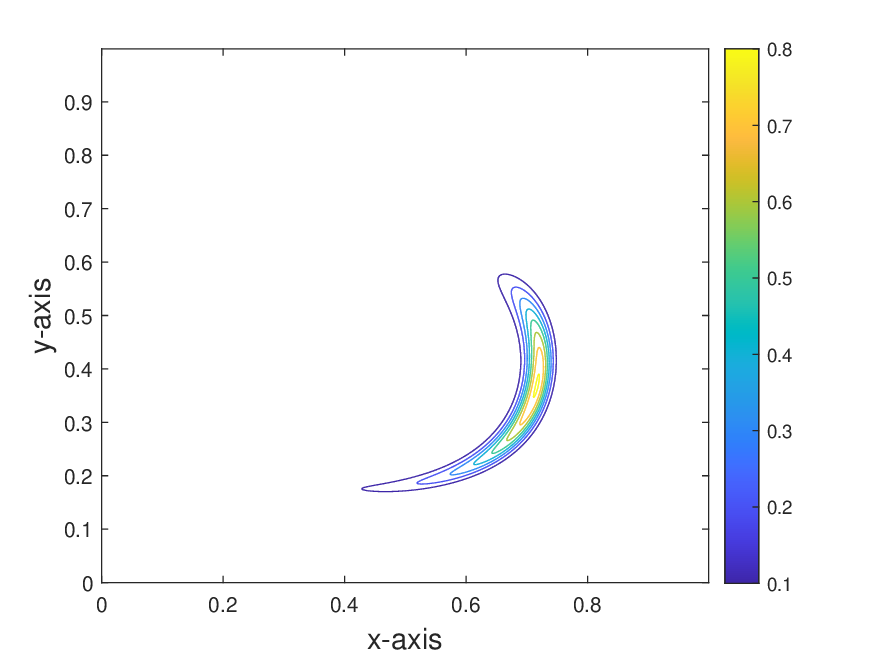}
	}
	\subfigure[$t=2$]
	{
		\includegraphics[width=0.3\textwidth]{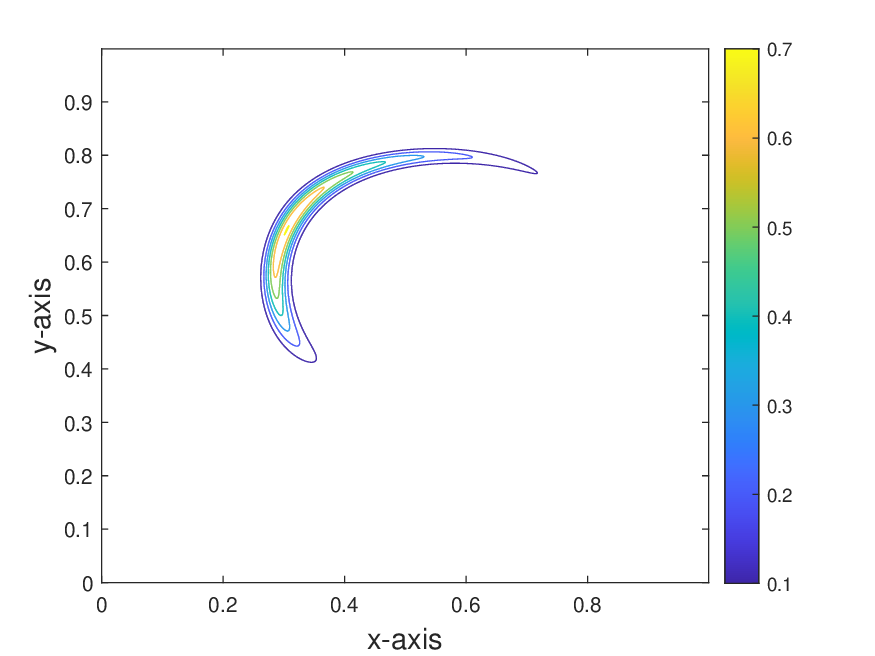}
	}
	\subfigure[$t=0.5$]
	{
		\includegraphics[width=0.3\textwidth]{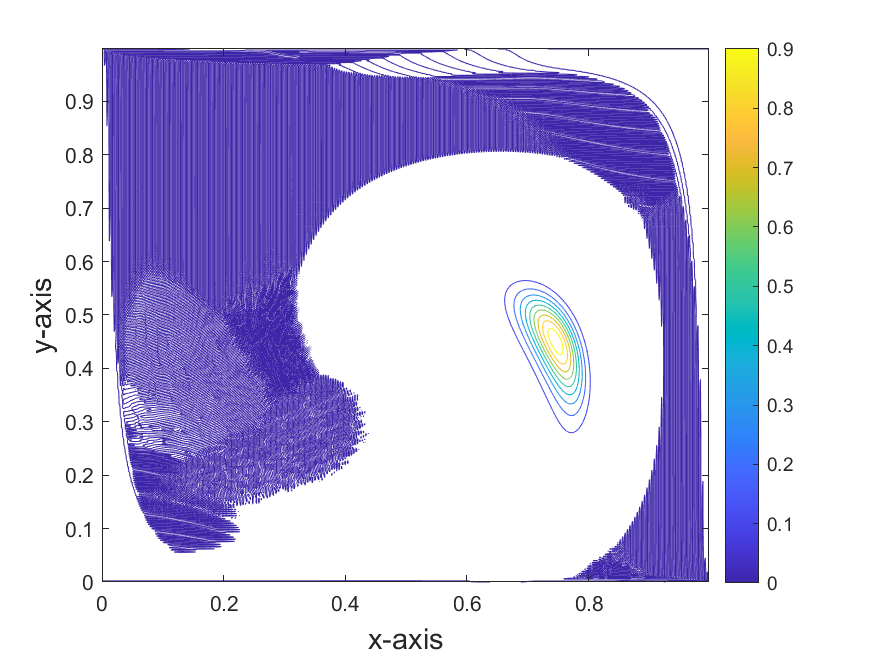}
	}
	\subfigure[$t=1.5$]
	{
		\includegraphics[width=0.3\textwidth]{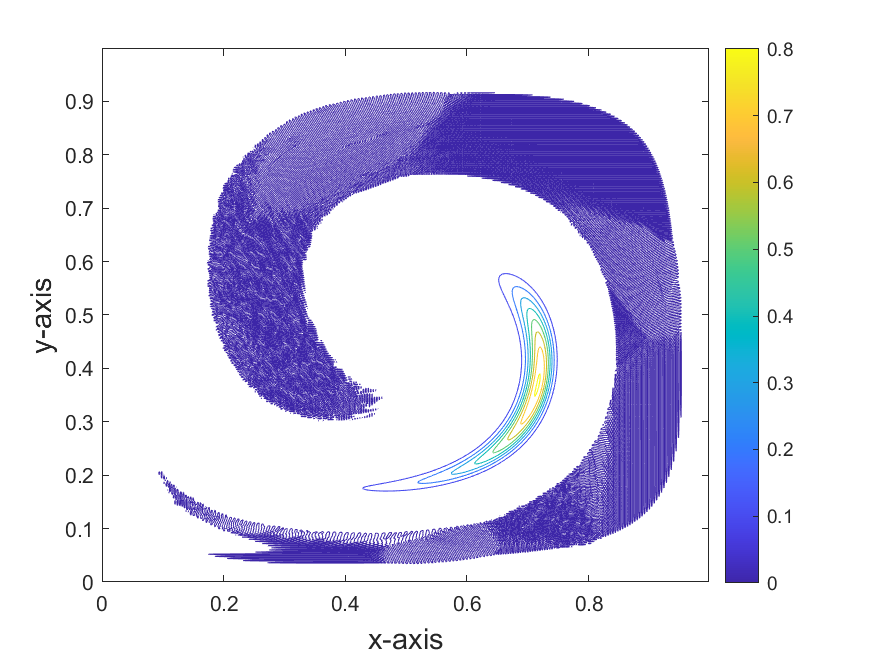}
	}
	\subfigure[$t=2$]
	{
		\includegraphics[width=0.3\textwidth]{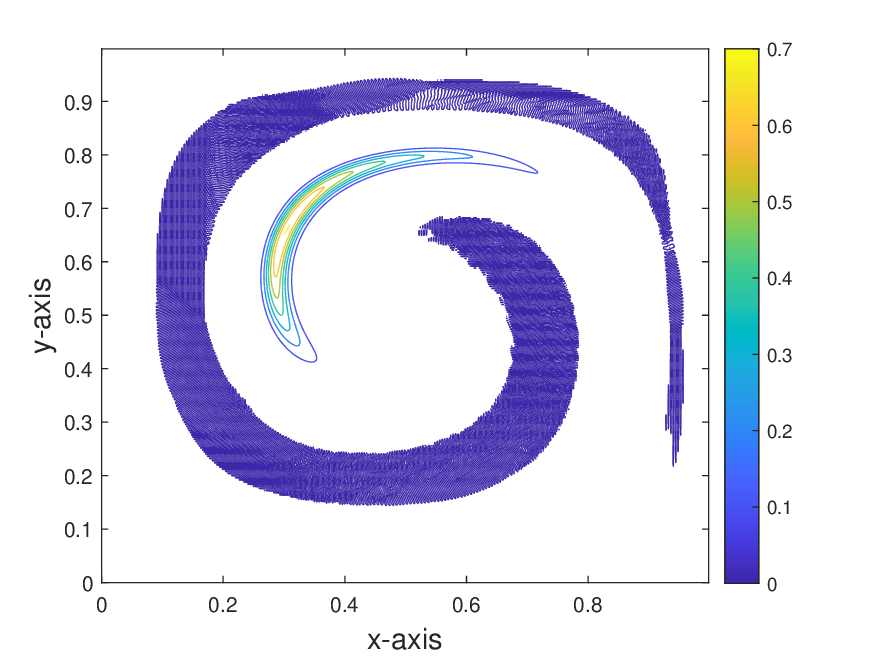}  
	}
	\subfigure[$t=0.5$]
	{
		\includegraphics[width=0.3\textwidth]{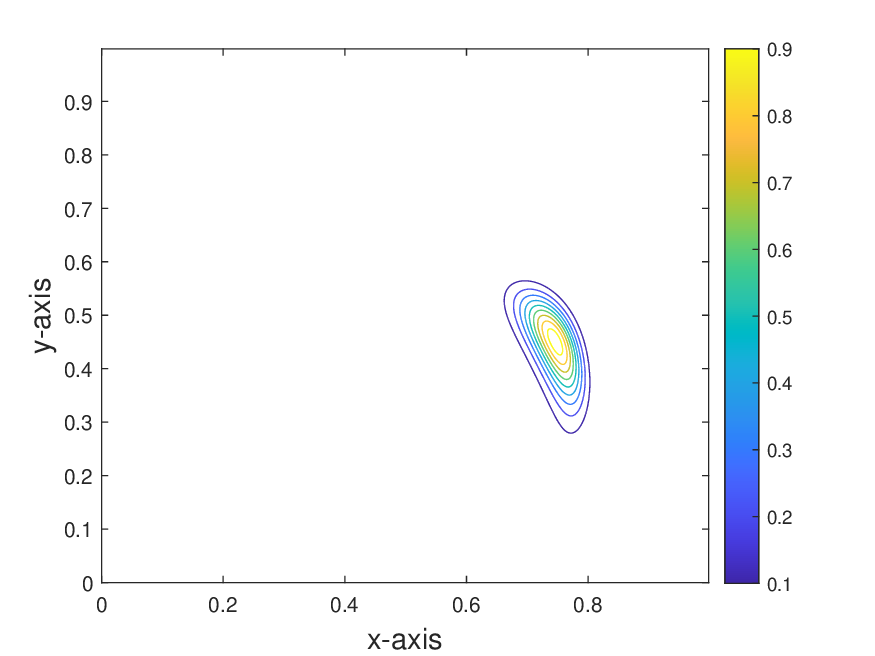}
	}
	\subfigure[$t=1.5$]
	{
		\includegraphics[width=0.3\textwidth]{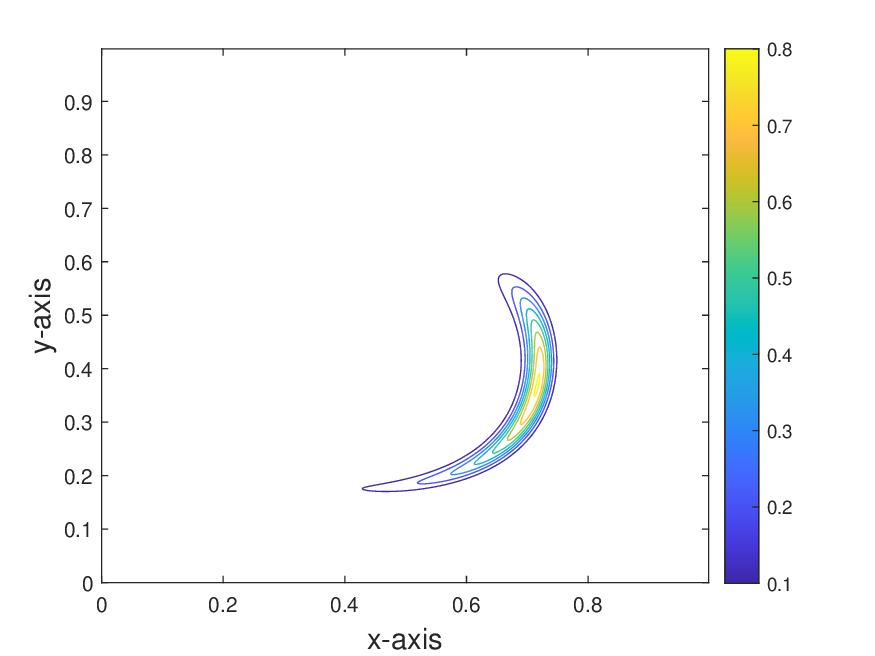}
	}
	\subfigure[$t=2$]
	{
		\includegraphics[width=0.3\textwidth]{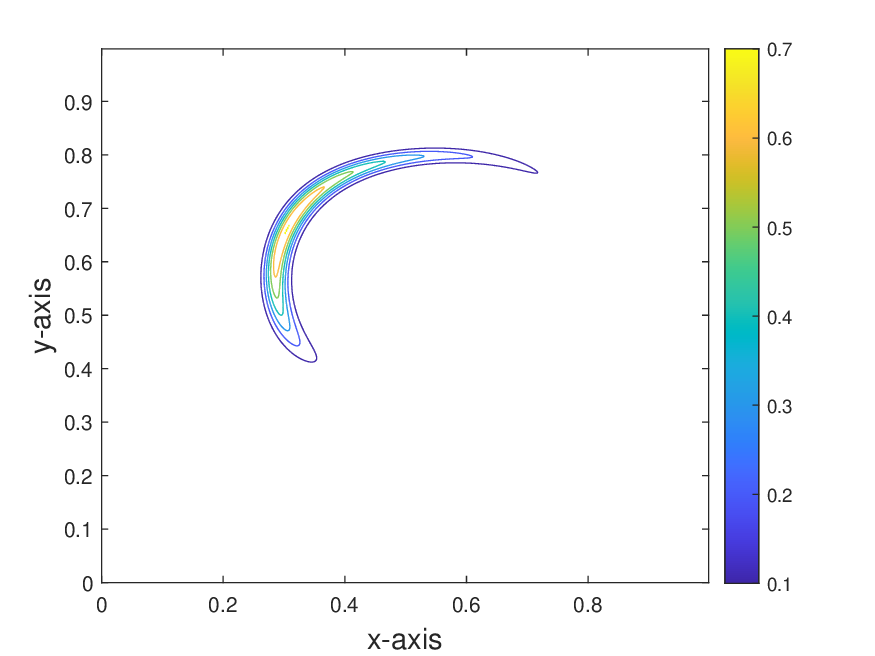}
	}
	\setlength{\abovecaptionskip}{0.05cm} 
	\setlength{\belowcaptionskip}{0.0cm}
	\caption{\small  Comparisons of contour plots of the reference solution (a-c), and the numerical solutions without the multipliers (d-f) and with multipliers (g-i)}
	\label{fig:roate:sqrt2:u}
\end{figure}
\begin{figure} [!htbp]
	\vskip -0.5cm

	\subfigure[{ HOC-ADI-Splitting }]
	{
		\includegraphics[width=0.43 \textwidth]{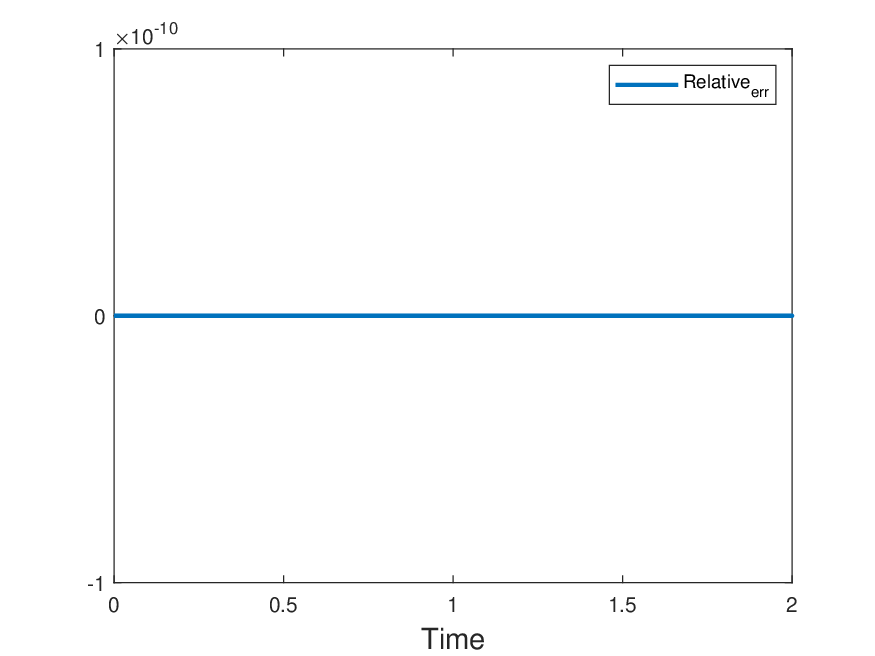}
	}
	\subfigure[ {BP-MC-HOC-ADI-Splitting}]
	{
		\includegraphics[width=0.43 \textwidth]{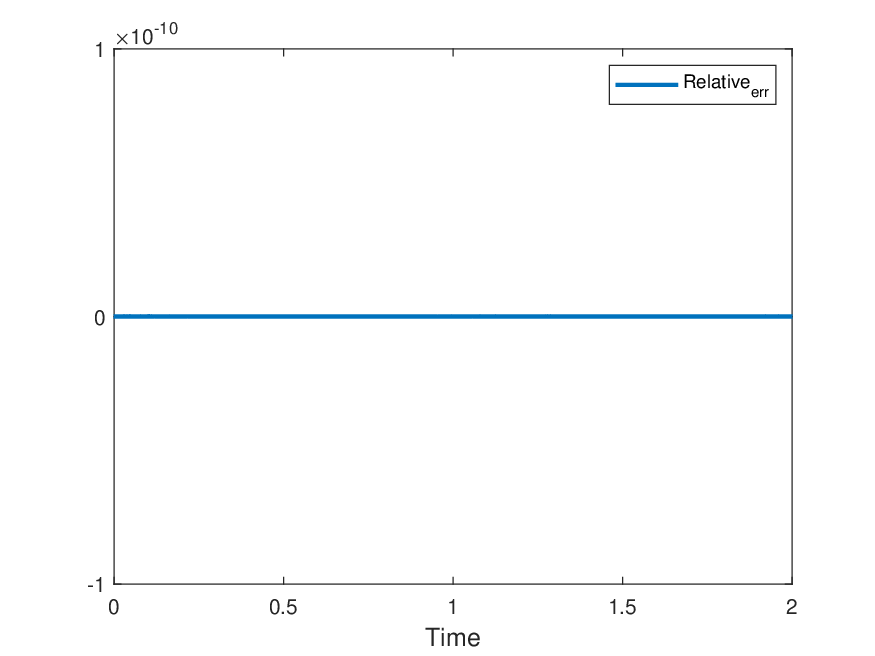}
	}
	\setlength{\abovecaptionskip}{0.05cm} 
	\setlength{\belowcaptionskip}{0.0cm}
	\caption{\small The numerical relative mass errors without/with Lagrange multipliers}
	\label{fig:roate:sqrt2:mass}
\end{figure}

In this simulation, the scalar quantity is gradually stretched and deformed into a spiral shape by the swirling flow when time marching. In Fig. \ref{fig:roate:sqrt2:u}, we show the contour plots of the numerical solutions obtained by the HOC-ADI-Splitting scheme and the BP-MC-HOC-ADI-Splitting scheme at $t=0.5$, $1.5$, $2$, respectively, by choosing $N=500$ and $\tau = \frac{h}{6{\max{|f^{\prime}|}}}$, where the reference solution is given by the BP-MC-HOC-ADI-Splitting scheme under $N=1000$ and $\tau = \frac{h}{6 { \max{|f^{\prime}|}}}$. As seen, the numerical solution obtained with the multipliers follows closely the reference solution, while that obtained without the multipliers is polluted with negative values and exhibits significant numerical dissipation, resulting in a less precise flow representation. Furthermore, the time evolution of the relative mass errors without/with the multipliers are plotted in Fig. \ref{fig:roate:sqrt2:mass}, which illustrates that the relative mass errors of the two proposed schemes approach the accuracy level of $10^{-10}$. Therefore, it again shows that discrete mass conservation is well maintained without/with the multipliers, which is consistent with the conclusion in Theorem \ref{thm:mass:2d}. 

\section{Concluding remarks} \label{sec_conclusion}
In this paper, we proposed several high-order compact and efficient splitting schemes for nonlinear convection diffusion equations. Generally speaking, our schemes leverage the advantages of the operator splitting method, the Lagrange multiplier approach, and the fourth-order compact difference approximations to address complex and multidimensional model problems with high-order accuracy. Apart from high-order accuracy and efficiency, the developed schemes allow us to achieve structure-preserving properties without complicated computations. Moreover, the proposed bound-preserving schemes automatically guarantee the uniform bounds of the numerical solution. In particular, for the 2D problems, high-order ADI splitting schemes are constructed, which further improve the computational efficiency. In addition, the substep method and the TVB limiter are also adopted to enhance the stability and effectiveness of the schemes. We also carried out a rigorous error analysis for the BP-HOC-ADI-Splitting scheme and derived optimal-order error estimates in discrete $L_2$ norm. It is hopeful that the error estimates could be generalized to the BP-MC-HOC-ADI-Splitting scheme, however, the proof is nontrivial and requires substantial new efforts. Numerical results demonstrate the accuracy, efficiency and structure-preserving properties of the proposed schemes.   
	
\section*{Declarations}
\begin{itemize} \small
	\item {\bf Funding}~ This work was supported in part by the National Natural Science Foundation of China (Nos. 12131014, 12371412, 11971482), by the Fundamental Research Funds for the Central Universities (No. 202264006) and by the  OUC Scientific Research Program for Young Talented Professionals.
	\item  {\bf Conflict of Interest}~ The authors declare that they have no competing interests.
	\item  {\bf Data Availability}~ Data will be made available on request.
\end{itemize}

\bibliographystyle{spmpsci}      
\bibliography{ref_L}

\end{document}